\documentclass[10pt]{article}
\usepackage{amsmath,amsfonts,amsthm,amssymb,amscd}

\catcode`\é=\active\def é{\'e} \catcode`\è=\active\def è{\`e}
\catcode`\ç=\active\def ç{\c{c}} \catcode`\à=\active\def à{\`a}
\catcode`\ê=\active\def ê{\^e} \catcode`\ù=\active\def ù{\`u}
\catcode`\ô=\active\def ô{\^o} \catcode`\û=\active\def û{\^u}
\catcode`\î=\active\def î{\^{\i}} \catcode`\â=\active\def â{\^a}
\catcode`\ö=\active\def ö{\"o}

\newcommand{\be}{\begin{equation}}
\newcommand{\ee}{\end{equation}}
\newcommand{\bs}{\begin{split}}
\newcommand{\es}{\end{split}}
\newcommand{\ba}{\begin{align}}
\newcommand{\ea}{\end{align}}
\newcommand{\basl}[1]{\begin{align}\begin{split}\label{#1}}
\newcommand{\bas}{\begin{align}\begin{split}}

\newtheorem{theo}{Theorem}[section]
\newtheorem{prop}[theo]{Proposition}

\newtheorem{defi}[theo]{Definition}

\def\carre{\hbox{\vrule \vbox to 7pt{\hrule width 6pt \vfill \hrule}\vrule }}

\newcommand\N{\mathbb{N}}

\newcommand\R{\mathbb{R}}
\newcommand\C{\mathbb{C}}

\title{Weyl calculus in QED I. The unitary group.}

\author{L. Amour\textsuperscript{1}, R. Lascar\textsuperscript{2} and
J. Nourrigat\textsuperscript{1}}

\date{ \ \textsuperscript{1}Université of Reims\hskip 1cm \ \textsuperscript{2}Université Paris 7, Denis Diderot  }

\begin{document}

\maketitle

\begin{abstract}
\noindent
In this work, we consider fixed $1/2$ spin particles interacting with the quantized radiation field in the context of quantum electrodynamics (QED). We investigate the time evolution operator in studying the reduced propagator (interaction picture). We first prove that this propagator belongs to the class of infinite dimensional Weyl pseudodifferential operators recently introduced in \cite {A-J-N} on Wiener spaces. We give a semiclassical expansion of the symbol of the reduced propagator up to any order with estimates on the remainder terms. Next, taking into account analyticity properties for the Weyl symbol of the reduced propagator, we derive estimates concerning transition probabilities between coherent states.
\end{abstract}

\parindent=0pt

\


\


\tableofcontents

\parindent = 0 cm

\parskip 10pt
\baselineskip 17pt

\section{Introduction.}\label{s1}

The aim of this article consists in studying a model of one or several fixed $1/2$ spin particles in a constant magnetic field ${\bf \beta}$ in the framework of quantum electrodynamics. This quantum model may be found in Reuse \cite{REU}. It is related to the spin-boson model, for which  spectral properties and some propagation results are known (see  \cite{H-S-1} \cite{A-H} \cite{A} \cite{H-S-2} \cite{DR-G-K}). A classical version of this model is introduced in \cite{B} in 1946. The main goal of this paper is to describe the operator $e^{-i{t \over h} H(h)}$ where $H(h)$ is the Hamiltonian of the model, using the infinite dimensional Weyl calculus
recently developed in
\cite{A-J-N} and \cite{A-L-N}. Potential benefits of this calculus is to perform analysis when the semiclassical parameter $h$ tends to $0$.

The Hilbert space describing the quantized field without interaction is usually defined as the symmetrized Fock space  ${\cal F}_s (H_{\C})$ over the complexified of a real Hilbert space $H$ (completion of the direct sum of symmetric tensor products, see \cite{RE-SI-1}). The space  $H_{\C}$ describes the states of the field with exactly one photon. The definition of this space is not chosen here as the usual one but rather inspired by the space given by
 Lieb-Loss \cite{L-L}. The definition is given in Section \ref{s2.A}. The
 Hilbert space of the system of $N$ spins ($1/2$) interacting with the quantized field is ${\cal F}_s (H_{\C})\otimes ({\C}^2)^{\otimes N}$.
The Hamiltonian
is recalled in Section \ref{s4}. Note that this operator is clearly simpler than the Pauli-Fierz operator where the particle is not
fixed. See, e.g., Cohen-Tannoudji, Dupont-Roc, Grynberg \cite{CT-DR-1}\cite{CT-DR-2}, Fr\"ohlich \cite{F-1}\cite{F-2}, Bach, Fr\"ohlich, Segal \cite{B-F-S} and Spohn \cite{SP}.

Our goal in this paper and in forthcoming works is to adapt to the operator $H(h)$ defined below in Section \ref{s4}, semiclassical results as those obtained in finite dimension, see e.g. \cite{Z}.
In order to do this, the Fock space
 ${\cal F}_s (H_{\C})$
is less convenient than an $L^2$ space to which it is isomorphic. The definition of this $L^2$ space is reminded in Section \ref{s2.B}. This space
is related to  an abstract Wiener space  $(i, H, B)$ where $H$ is the real Hilbert space, $B$ is a suitable Banach space containing $H$ and $i$ the injection of
$H$ into $B$. The precise assumptions fulfilled by $B$ are recalled in Section \ref{s2.B}. For every
$h>0$, the space $B$ is endowed with a gaussian measure with variance $h$ denoted by $\mu _{B , h}$.
The main properties are recalled in \cite{A-J-N}. The space ${\cal F}_s (H_{\C})$ is then isomorphic
to $L^2 (B , \mu _{B , h/2})$ (see \cite{SI} \cite{J})) also denoted here  ${\cal H}_{ph}$. The choice of $h/2$ as the variance is taken in order to agree with usual formulas of pseudodifferential calculus. Taking account of this isomorphism, the Hilbert space of the system considered here becomes $L^2 (B , \mu _{B , h/2})\otimes ({\C}^2)^{\otimes N} = {\cal H}_{ph} \otimes ({\C}^2)^{\otimes N}$. The two spaces $H$ and $H_{\C}$ are specified in Section \ref{s2.A}. The space $H_{\C}$ (or $H\times H$)  is the single photon space.

Two types of operators act in this space. The Hamiltonian  $H_{ph}$ of the (non interacting) free field
is an unbounded selfadjoint operator in ${\cal F}_s (H_{\C})$ or in $L^2 (B , \mu _{B , h/2})$, not described with the infinite dimensional Weyl calculus. It is usually defined as $H_{ph} = {\rm d}\Gamma (h M_{\omega})$ where the functor $\Gamma $ and ${\rm d}\Gamma $ are standard operators defined in Fock spaces
 (see \cite{RE-SI-2}) and $M_{\omega}$ is acting in the one photon space $H_{\C}$  (see Section \ref{s3.B}). The notation $H_0 = H_{ph} \otimes I$ is also used here. The operator
$e^{i{t\over h}  H_{ph}} $ equals $\Gamma ( e^{it M_{\omega}})$.
This operator may be also viewed as a kind of metaplectic transform associated with the symplectic linear transform  $\chi_t$ in $H^2$, which is also unitary
(see Section \ref{s3.B} and see \cite{LA} for a metaplectic group in this context).

The other operators involved here are pseudodifferential operators. Weyl calculus associates bounded and unbounded operators in
 ${\cal H}_{ph} = L^2 (B , \mu _{B , h/2})$ to suitable functions (symbols) defined on  $H^2$. The definition of the operator associated with a symbol $F$ being a continuous linear form on $H^2$ is standard: it is usually called {\it Segal field}
and defined with creation and annihilation operators. The operators $B_j(x)$ and $E_j(x)$ ($1\leq j \leq 3$) associated with the components of the magnetic and electric fields at the point $x$ are defined in this way, for all  $x$ in $\R^3$. See \cite{F-1}\cite{F-2}\cite{SP} and \cite{RE-SI-2}.

The Weyl calculus developed in \cite{A-J-N} or \cite{A-L-N} is an extension of the notion of Segal field. It enables to associate an operator (bounded or not) $OP_h^{weyl} (F)$ with some functions $F$ on $H^2$, for each $h>0$. To define a bounded operator, it is sufficient that the function  $F$ belongs to a class $S_2({\cal B}, \varepsilon )$ introduced in Definition \ref{d3.1} below, associated with an Hilbertian basis ${\cal B} = (e_j)$ of $H$ and with a summable sequence $\varepsilon = (\varepsilon_j)$.  If the symbol
$F$ takes values in ${\cal L} ( ({\C}^2)^{\otimes N})$ then the corresponding operator acts in ${\cal H}_{ph}\otimes ({\C}^2)^{\otimes N}$.
Theorems \ref{t1.1} and \ref{t1.2} below essentially rely
 on Beals characterization Theorem (Theorem 1.2 of \cite{A-L-N}).

The operator $H(h)$ denotes the Hamiltonian of the entire system, for all $h>0$ and  $H_0 = H_{ph} \otimes I$.
These two operators are selfadjoint and share the same domain. The operator $H(h)$
is precisely defined in Section \ref{s4}.

The first main result is the following one.

\begin{theo}\label{t1.1} Let $H(h)$ be the Hamiltonian defined in Section \ref{s4} and $H_0 = H_{ph} \otimes I$.
Then, there exists a function $U(t,q,p,h)$ defined on $\R \times H^2$ with parameter $h>0$ and taking  matrices values (in
${\cal L} ( ({\C}^2)^{\otimes N})$) fulfilling,
\be\label{1.1}
e^{-i{t \over h}  H(h) } = e^{-i{t \over h}   H_0} U_h^{red} (t),\qquad U_h^{red} (t) =
Op_h^{weyl} (U(t, \cdot , h)) \ee
%---
and having the following properties.

$(i)$ As a function of $(q , p)$, it belongs to the space
$S_{\infty}^{mat} ({\cal B}, |t| \varepsilon (t))$ introduced in Definition \ref{d3.1} associated with the Hilbertian basis ${\cal B}$ of Section \ref{s2.C} and the sequence $\varepsilon (t) $ of Proposition \ref{p5.1}.

$(ii)$ Additionaly, this function is bounded in this space (in the sense of Definition \ref{d3.1}) as $h$
runs over $(0, 1)$ and $t$ belongs to a compact set of $\R$.

$(iii)$ Moreover,
the function $U ( t, q, p, h)$ has a asymptotic expansion in powers of  $h$ up to any order $m$,
%---
\be\label{1.1b} U ( t, q, p, h) = \sum _{j= 0 }^m g_j ( t, q, p) h^j + h^{m+1} R_m (t, q, p, h), \ee
%----
where the $g_ j (t, \cdot)$ and $R_m(t, \cdot , h)$ belong to $S_{\infty}^{mat} ({\cal B}, K|t| \varepsilon (t))$, for some non negative real constant $K$. These functions are solutions to the differential system (\ref{6.1a})(\ref{6.1b}).

$(iv)$ In addition, the function $(q, p) \mapsto U(t, q, p, h)$
is analytic on $H^2$ and and there exits an holomorphic extension on $H_{\C}^2$ taking values in ${\cal L} (( \C^2  )^{\otimes N} )$ and denoted ${\cal U} (t, q, p, h) $ which satisfies,
%----
\be\label{1.2} |   {\cal U} (t, q, p, h) |  \leq M(h, t) e^{ K |t| |{\rm Im } (q , p) |},\qquad (q , p) \in H_{\C}^2 \ee
%---
where $M(t, h)>0$ is bounded when $h$ belongs to $(0, 1)$ and $t$ to a compact in $\R$, and  $K$ is some non negative real constant.
\end{theo}

The proof uses a Calder\'on-Vaillancourt type result proved in \cite{A-J-N} and a Beals characterization type theorem obtained in \cite{A-L-N}.

Note that the parameter $h$ appearing in the Hamiltonian $H(h)$ is also the parameter related to the variance of the Gaussian measure $\mu_{B,h/2}$ defined above and associated to the $L^2$ space.

In the case of a single spin ($N=1$), the system (\ref{6.1a}) satisfied by  $g_0$ is related to Larmor precession.  The expansion (\ref{1.1b}) is different from perturbative series (see, e.g., Section 8.3 of \cite{G-R}, formula (4.175) in \cite{REU}, Section 4.2 in \cite{PS}, Section 3.5 of \cite{W}, $\dots$), expansions in powers of the coupling constant. Indeed, the expansion (\ref{1.1b}) is a semiclassical expansion whereas perturbative series are not.

In finite dimension, the semiclassical Weyl calculus enables
for instance to control the remainder terms up to any order of semiclassical expansions of various quantities or functions related to the Hamiltonian under consideration (for example, see \cite{C-R}\cite{D-S}\cite{H}\cite{M}\cite{R}\cite{Z}). In (\ref{1.1b}), the precise control of the remainder terms relies on the estimates on the derivatives of $U(t, q, p, h)$ which are in Definition \ref{d3.1} of $S_{\infty}^{mat} ({\cal B}, |t| \varepsilon (t))$ (See Section \ref{s6}).
Recently, the definition and some properties of the Weyl calculus are extended to the Wiener spaces setting: norm estimates (\cite{A-J-N}),  Beals type characterization (\cite{A-L-N}). These properties are involved in the present work.
One may think that others properties will be soon also extended to the Wiener spaces setting and the resulting applications may concerned the semiclasical expansion and control of the remainder terms of the time evolution for some average values of some physical observables (Ehrenfest type result). This study should probably be local in time, as in finite dimension. Another tool for semiclassical analysis comes from coherent states depending on the parameter $h$ (see definition in Section \ref{s3.A}) and their semiclassical evolution for the Hamiltonian $H(h)$ defined in Section \ref{s4} is given in \cite{A-N}.  Here, the applications are only the expansion (\ref{1.1b}) and the  result below (Theorem \ref{t1.2}) which is a consequence of the analyticity properties of the symbol.

Theorem \ref{t1.1} $(i)(ii)(iv)$ is proved in Section \ref{s5} and the semiclassical expansion of the symbol $U(t, q, p, h)$ (function  defined on $\R \times H^2$ and taking values
 in ${\cal L} ( ({\C}^2) ^{\otimes N})$) of Theorem \ref{t1.1} $(iii)$ is given in
Section \ref{s6}.

Operators with analytic symbols are introduced in \cite{BM-K} in finite dimension.
An application of the analyticity properties is to derive estimates of transition probabilities between coherent states.
Coherent states $\Psi_{X , h}$ are elements of the Fock space, or alternatively of the $L^2$ space which is isomorphic, indexed by an element $X = (a , b)$ of
 $H^2$ where $H$ is the one photon real Hilbert space. Moreover, the coherent states depend on the semiclassical parameter $h>0$ and their definitions are reminded in (\ref{csf})(\ref{csl}).
In the case of finite dimension, if a propagator is known to be the composition of a
  pseudodifferential operator with a metaplectic operator as in  (\ref{1.1}) then an estimate of the transition probability between coherent states is given in \cite{U}. The analyticity of the symbol here improves this estimate.  Let us recall that the symplectic map $\chi_t$ in $H^2$  defined in (\ref{3.14}) reflects the free dynamics of photons between times $0$ and $t$.

Theorem \ref{t1.2} below is proved in Section \ref{s7}.

\begin{theo}\label{t1.2}  For all $X$ and  $Y$ in $H^2$, for any $a$ and $b$
 in  $( {\C}^2 )^{\otimes N}$ with norm $1$, the following estimate holds,
\be\label{1.3}
 | < e^{i{t \over h}  H(h)} ( \Psi _{X h}\otimes a )  , \Psi _{Y h}\otimes b > | \leq M( h, t)
  e^{  K |t| |X-\chi_t^{-1} (Y)|- {1\over 4h} |X- \chi_t^{-1} (Y) |^2 }
  \ee
%---------
where $\chi_t$ is the symplectic map in $H^2$ defined in (\ref{3.14}) and where
$M( h, t)$ and $K$ are the constants of Theorem \ref{t1.1} $(iv)$.
\end{theo}

The benefit of Theorem \ref{t1.2} comes from the comparison with standard results. It is proved in Unterberger \cite{U} that, if $A_h$ is an operator associated by the Weyl calculus in finite dimension with a function which is bounded together with all of its derivatives and if the $\Psi_{X h}$ ($X \in \R^{2n}$)  are the usual coherent states then, for all integers $N$, there exists $C_N>0$ such that, for all $X$ and $Y$ in $\R^{2n}$,
 for any $h$ in $(0, 1)$,
 %-----
$$ |<A \Psi_{X h}, \Psi_{Y h} > | \leq  C_N  \left [ 1 + {|X-Y| \over \sqrt h} \right ] ^{-N}. $$
%--
Note the two differences with Theorem \ref{t1.2}. Namely, the free photons evolution leads to replace $Y$ by $\chi_t (Y)$ and the analytic properties of the symbol enables to improve the rapid decay estimate which is now replaced by an exponential decay estimate.

A model for a fixed  $1/2$ spin particle in a constant magnetic field interacting with the quantized field and also interacting with a rotating (non quantized) magnetic field may also be  considered. It is seen in Section \ref{s4} that such a model may be reduced with a suitable transform to the case where only the constant field is not quantized. Thus, a standard method in quantum mechanics remains valid in QED (see Theorem \ref{t4.2}).

{\it Notations.} The notation $( \ , \ )$ stands here for real scalar products on $H$ or $H^2$ and $< \ , \ >$ denotes Hermitian products, antilinear with respect to the second variables. These products are used on  ${\cal H }_{ph}$, $({\C}^2)^{\otimes N}$,
and on $H^2$ when identified to the complexified $H_{\C}$. In that case, it is written as,
\be\label{1.4}
< X , Y> =   (X , Y)  + i \sigma (X , Y ),\qquad X, Y \in H_{\C}.\ee
%---

\section{Usual Hilbert spaces in  QED.}\label{s2}

\subsection{The one photon Hilbert space.}\label{s2.A}

The space $H$ (resp. $H_{\C}$) stands for the set of mappings
 $f = (f_1, f_2, f_3)$
with $f_j$ belonging to $L^2({\R}^3)$ and taking real (resp. complex) values, satisfying
\be\label{2.1} k_1f_1(k) + k_2f_2(k) + k_3f_3(k) = 0.\ee
This space is endowed with the norm
\be \label{2.2} |f|^2 = \sum_{j=1}^3 \int_{\R^3} |f_j(k)|^2dk.\ee
The space $H_{\C}$ is here the Hilbert space corresponding to a unique photon. It is the space of divergence free for Fourier transforms of vector fields. As later seen, it is also related to the set of all possible initial data for Maxwell equations in vacuum (without interaction).

Let us also mention a more common way to introduce this space. For all $k\in\R^3$ such that $k_1^2 + k_2^2 \not= 0$, one can find $e_1(k)$ and
$e_2(k)$ in $k^\bot$, being an orthonormal basis of $k^\bot$ and $C^\infty$ regular with respect to $k$ in its domain of definition. For each $f$ in
$H$ and for any $k$ in its domain of definition, set $f_\lambda(k) =
f(k) \cdot e_\lambda(k)$ $( \lambda= 1, 2)$. One has,
$$f(k) = f_1(k)e_1(k) + f_2(k)e_2(k)$$
and
$$ |f|^2 = \sum_{\lambda=1,2}\int_{R^3} |f_\lambda(k)|^2dk,$$
for any $k$ in the domain.
Thus, $H$ is isomorphic to $L^2({\R}^3)^2$ or to $L^2({\R}^3\times\{1,2\})$. This is the usual way but less convenient, since a singularity now appears in the set of  $k\in\R^3$ with $k_1^2+k_2^2=0$.

\subsection{The quantized field Hilbert space.}\label{s2.B}

It is the symmetrized Fock space ${\cal F}_s(H_{\C})$ over $H_{\C}$. For  pseudodifferential analysis, we rather use an $L^2$ space, to which it is isomorphic. To this end, a Banach space $B$ containing
$H$  may be found, together with a suitable measure.
Let us first recall the hypotheses fulfilled by  $B$. For all finite dimension subspaces $E$ in $H$ and for all $h>0$, a measure $\mu _{E , h}$ is defined by setting,
%------
\be\label{2.3} \mu _{E , h} (  \Omega)  = (2\pi h)^{-{\rm dim} (E)/2}
 \int _{\Omega}   e^{-{|y|^2 \over 2h}} dy,\ee
%----
for all Borel sets $\Omega $ in $E$.

\begin{defi}(\label{d2.1}\cite{G-2}\cite{G-3}\cite{K}) Let $H$ be real separable Hilbert space with the norm $|\cdot|$. Define on  $H$  another norm
$N$ satisfying for some constant  $C > 0$,
%---
\be\label{2.4} N(x) \leq  C |x|.\ee
%-----
The norm $N$ is called measurable if for all $\varepsilon > 0$ and for all $h > 0$, there exists a finite dimensional subspace $F$ in $H$ such that, for any finite dimensional subspaces $E$ orthogonal to $F$,
\be\label{2.5} \mu_{E,h}(\{x\in E,\ N(x)>\varepsilon\}) < \varepsilon  \ee
%---
where  $\mu_{E,h}$ is defined in (\ref{2.3}).
\end{defi}

Such a norm can always be found. The completion $B$ of $H$ for such a norm satisfies,
%---
\be\label{2.6} B' \subset H \subset B,\ee
%---
each space being dense in the next one.  Denoting $i$ the injection from $H$ into $B$,
 $(i, H, B)$ is called a Wiener space. Then, for all $h>0$, the Borel $\sigma-$algebra on $B$ is equipped with a measure denoted by $\mu _{B , h}$, related through natural formulas to the measures $\mu _{E , h}$ defined in (\ref{2.3}), for all finite dimensional subspaces $E$ of $H$.
 Then, the Fock space ${\cal F}_s(H_{\C})$ becomes isomorphic to $L^2(B,\mu_{B,h})$, for all $h$.
 To fit with standard formulas of pseudodifferential calculus, we choose
 $L^2(B,\mu_{B,h/2})$ as the Hilbert space of the quantized field, isomorphic to the Fock space.

We now specify the choice of space $B$ suitable for the space $H$ defined in Section \ref{s2.A}. It is already known that,  if $D$ is
an injective Hilbert-Schmidt operator in a Hilbert space $H$, then the
norm $x \rightarrow |Dx|$ is measurable on $H$ (\cite{G-3}), example 2, p. 92).
In our case, with $H^0$ being the space $L^2({\R}^3)^3$, the operator
%----
\be\label{2.7} D = (I - \Delta_k + |k|^2)^{-m},\qquad  m > {3 \over 2} \ee
%----
acting on the three components, is Hilbert Schmidt in $H^0$, and the norm
$f\mapsto N(f)=|Df|$ is therefore measurable on $H^0$.
It is easily seen that a measurable norm restricted to a closed subset is a measurable norm on that subset. The subset $H$ of $H^0$ is closed. Consequently, one may choose as a space $B$, the completion of  $H$ defined above for the norm $N$. Thus $(i, H ,
B)$ is a Wiener space.

Usually in field theory (see \cite{SI}), the space ${\cal S}'(\R^3)$ is equipped with a probability measure and the space $L^2 ( {\cal S}'(\R^3) )$ is isomorphic to  ${\cal F }_s ( H_{\C} ) $. In order to use the results in \cite{A-J-N} and \cite{A-L-N}, which themselves use results of  Gross, Kuo and Janson (\cite{G-1}\cite{G-2}\cite{G-3}\cite{J}\cite{K}), the space ${\cal S}'(\R^3)$ has to be replaced with a suitable Banach space $B$
 containing the space $H$ of Section \ref{s2.A}. It is on the Borel $\sigma-$algebra of this space $B$ (which is not unique) that is defined the probability measure that is used here. If instead of the space defined in Section \ref{s2.A}, we simply had $L^2(\R^3)$, then the space $B$ would be a negative order Sobolev space sufficiently large, actually relatively close to ${\cal S}'(\R^3)$.

For all  $a$ in $H$,
one uses a function $\ell _a $ defined on  $B$ in the following way.
When $a$ belongs to $B' \subset H$, one has $\ell_a(x) = a(x)$. If $a$ is $H$, it is approximated by a sequence $(a_j)$ in $B'$. The sequence $\ell_{a_j}$ is proved to be a Cauchy sequence in
 $ L^2(B , \mu_{B , h/2})$, and we denote by $\ell _a$ its limit.

The following subspace ${\cal D}$ of $L^2(B,\mu_{B,h/2})$ is often used in the sequel.

\begin{defi}\label{d2.2} For all finite dimensional subspaces $E$ of $H$, ${\cal D}_E$ denotes the space of functions $f : B \rightarrow
\C$ satisfying,

i) The function $f$ is written as $g \circ P_E$ where $g$ is a continuous function from $E$ to $\C$ and where $P_E$ is the mapping from $B$
to $E$ defined as following, when choosing an orthonormal basis $\{u_1,\dots,u_n\}$ of $E$,
%---
\be\label{2.8} P_E(x) = \sum_{j=1}^n \ell_{u_j}(x) u_j,\quad {\rm for\ a.e.\ } x\in B\ee
%----
(This mapping $P_E$ is independent on the chosen basis).

ii) The function $E^2 \ni X \mapsto < f, \psi_{X h} >$ belongs to the Schwartz space
${\cal S}(E^2)$.

The union of all subspaces ${\cal D}_E$ is denoted by {\cal D}.
\end{defi}

According to Definition 4.4 of \cite{A-J-N}, we say that a continuous function $f$ on $H$
has a stochastic extension denoted $\widetilde f$ in $L^p (B , \mu _{B , h})$  ($1 \leq p < \infty $) if,
for all non decreasing sequence of finite dimensional subspaces $(E_n)$ in $H$ with a dense union in $H$, the sequence $f \circ P_{E_n}$  is in $L^p (B , \mu _{B , h})$
and converges to $\widetilde f$ in this space. Examples are found in Section 8.2 of \cite{A-J-N}.

The next property is used in the following.

\begin{prop}\label{p2.3} Let $f : H \rightarrow {\C}$ be a continuous function with a stochastic extension  $\widetilde f $ in $L^1(B , \mu _{B , h})$. Let
$\chi : H \rightarrow H$ be a unitary continuous linear mapping. Then
$g = f \circ \chi $ has a stochastic extension $\widetilde g$ in
$L^1(B , \mu _{B , h})$ and one has,
%---
\be\label{2.9}\int _B \widetilde g (y) d\mu _{ B , h} (y) = \int _B \widetilde f (x)
d\mu _{ B , h} (x).\ee
%----
\end{prop}

{\it Proof.} Let $(E_n)$ be a non decreasing sequence of finite dimensional subspaces in $H$ whose union in dense in $H$. Set $F_n = \chi  (E_n)$. The sequence
$(F_n)  $ shares the same properties.  The hypothesis on $f$ implies that the sequence
$ f \circ P_{F_n }$ is a Cauchy sequence in $L^1 ( B , \mu_{B , h})$. For all finite dimensional subsets $X$ and $Y$
of $H$ with $X \subset Y$, let $\pi _{XY}$ be the orthogonal projection from $Y$ to $X$.  According to the  transfer Theorem, the assumption gives  when $m < n$,  that the integral,
%---
$$ I_{mn} = \int _{F_n } | f(x) - f ( \pi _{F_m, F_n} (x)) |  d\mu_{F_n , h} (x) $$
%----
goes to 0 as $m =\inf ( m, n)$ tends to infinity. The change of variables $x = \chi (y) $ shows that,
%---
$$ I_{mn} = \int _{E_n } | g (y) - g  ( \pi _{E_m, E_n} (y)) |  d\mu_{E_n , h} (y). $$
%----
Using the  transfer Theorem, this proves that $ g \circ P_{E_n }$ is a Cauchy sequence in
$L^1 ( B , \mu_{B , h})$. Since the sequence $E_n$ is arbitrary then the function $g$ has a stochastic
extension in $L^1 ( B , \mu_{B , h})$. For all $n$, one has,
%----
$$
\int_{B^2} (G\circ\chi)(\tilde{\pi}_{E_n}(Z))d\mu_{B^2,h/2}(Z) =
\int_{E_n} (G\circ\chi)(Z)d\mu_{E_n,h/2}(Z)$$
%----
$$
=
\int_{\chi E_n} G(Z)d\mu_{\chi E_n,h/2}(Z)
=
\int_{B^2} G(\tilde{\pi}_{\chi E_n}(Z))d\mu_{B^2,h/2}(Z).$$
%------
We deduces (\ref{2.9}) making $n$ going to infinity.\hfill$\Box$

\subsection{An Hilbertian basis of $H$.}\label{s2.C}

The set of pseudodifferential operators recalled in Section \ref{s3.A} is not invariant by a change of basis. We then have to choose a particular basis of $H$. Here $E$ denotes the set of $L^2$ vector fields on the unit sphere $S^2$, namely, the mappings $f = (f_1, f_2, f_3)$ on $S^2$, taking values in  $\R^3$ and satisfying,
$$\omega_1 f_1 (\omega) + \omega_2 f_2 (\omega) +\omega_3 f_3 (\omega)  = 0 {\rm \ a.\ e.}$$
This space is equipped with the norm
$$||f||^2_E = \sum_{j=1}^3 \int_{S^2} |f_j(\omega)|^2 d\mu(\omega)$$
where $\mu$ is natural measure on $S^2$ ($\mu(S^2)=4\pi$). Thus $H$
is identified to $L^2(\R_+, r^2dr, E)$. We obtain an Hilbertian basis of $H$ written as,
%-----
\be\label{2.10} f_{mn} (k) = u_m(|k|) v_n(k/|k|)\ee
%----------------
where $(u_m)$ is an Hilbertian basis of $L^2(\R_+, r^2 dr)$ and  where $(v_n)$ is an Hilbertian basis of  $E$.
We are then led to make a particular choice.

The space $E$ may be identified as the space of  1-differential forms on the unit sphere $S^2$, being in  $L^2$ for the natural measure on this sphere. The de Rham Laplacian $\Delta _S$ is available on this space. One chooses the basis $v_n$ of $E$ as an Hilbertian basis of $E$ constituted with of eigenvectors of $\Delta _S$,
%---
$$ \Delta _S v_n = \mu_n v_n.$$
%------
The $\mu_n$'s are the non decreasing sequence of eigenvalues of the de Rham Laplacian, repeated with their multiplicity. The function $u_m$ are the eigenvectors of the following operator,
%-----
$$ L = - {d^2 \over dr^2} - {2\over r} {d \over dr} + r^2 \leqno(2.11)$$
%------
being essentially selfadjoint in  $L^2(\R_+, r^2 dr)$ (see Reed-Simon \cite{RE-SI-2}, Theorem
X.11). One checks that,
$$ L u(r) = \lambda u(r), $$
 setting
$u(r) = e^{-{r^2 \over 2}} \psi(r^2)$ and $x = r^2$, is equivalent to,
%----
\be\label{2.12}x\psi''(x) = \Big( {3\over 2} - x\Big)\psi'(x) + {\lambda-3 \over 4} \psi(x) = 0.\ee
%-------
Thus, the eigenvalues of $L$ are  the $\lambda_m = 4m+ 3$ ($m \in\N$) of multiplicity $1$ and the corresponding functions $\psi_m$ in $(\ref{2.12})$ are the generalized Laguerre's polynomials
 $L_m^{(1/2)}$,
up to a multiplicative factor. There exists an Hilbertian basis $(u_m)$ of $L^2(\R_+, r^2 dr)$
constituted with eigenfunctions of $L$,
 %----
\be\label{2.13} L u_m = (3 + 4m)u_m,\qquad u_m(r) = C_m e^{-{r^2 \over 2}} L_m^{(1/2)}(r^2).\ee
%-----

When the Banach space $B$ is constructed as before, one sees that the $f_{mn}$ are elements of $B'$. Indeed, the three components of $f_{mn}$ are in ${\cal S}(\R^3)$ and consequently, for all $\varphi$ in $H$,
$$ | ( f_{mp},\varphi ) | \leq C_{mp}|D\varphi|$$
where $D$ is defined in (\ref{2.7}).

\begin{prop}\label{p2.4}  For all mappings $\varphi = (\varphi_1 , \varphi_2 , \varphi_3)$ in
$H$ with components in ${\cal S} (\R^3)$ and vanishing in a neighborhood of the origin, one has,
%---
$$ m^{\alpha } n^{\beta }  |(f_{mn}\cdot \varphi) | \leq  \sum _{j=1}^3 N_{\alpha \beta } (\varphi _j)
$$
%---
where $N_{\alpha \beta }$ is a semi-norm on ${\cal S}(\R^3)$. In particular, the sequence $|(f_{mn})\cdot \varphi |$
is summable.
\end{prop}

{\it Proof.} One has,
%---
$$ \lambda _m^{\alpha } \mu_n ^{\beta} (f_{mn})\cdot \varphi = ( L^{\alpha } u_m
\otimes \Delta_S ^{\beta } v_n ) \cdot \varphi  = f_{mn} \cdot (( L^{\alpha } \otimes \Delta _S ^{\beta }) \varphi ).$$
%----
If $\varphi $ belongs to ${\cal S} (\R^3 , \R^3)$ and vanishes in a neighborhood of the origin, then integrating by parts gives,
%----
$$ \lambda _m^{\alpha } \mu_n ^{\beta}  | (f_{mn})\cdot \varphi | \leq \Vert
(( L^{\alpha } \otimes \Delta _S ^{\beta }) \varphi ) \Vert. $$
%---
Since $ m \leq \lambda _m = 4m + 3$ and since the eigenvalues  $\mu_n$ of the de Rham Laplacian
 satisfies $n \leq C \mu_n$ for some constant $C$  (see Ikeda-Tanigushi \cite{I-T} or Folland \cite{F}) then the proof of
the Proposition is completed.

\hfill $\Box$

\section{Usual operators in QED.}\label{s3}

\subsection{Pseudodifferentials operators.}\label{s3.A}

\begin{defi}\label{d3.1} Let ${\cal B} = ( f_{jk})$ be the Hilbertian basis of $H$
defined in Section \ref{s2.C}.
We called  multi-index  $(\alpha, \beta )$ a map from ${\N}^2$ into
 $\N \times \N$ satisfying $\alpha_{jk} = \beta_{jk}= 0$ for all but a finite number of indices.
For all integers $m\geq 0$,  ${\cal M}_m$  denotes the set of multi-indices $(\alpha , \beta)$ satisfying $\alpha _{jk} \leq m$ and
$\beta _{jk} \leq m$, for all $(j, k)$.
For any integer $m\geq 0$, any $M>0$, and for each summable sequence of nonnegative real numbers  $\varepsilon = (\varepsilon_{jk} )$, we denote by $S_m ( {\cal B} ,M,  \varepsilon )$ the set of continuous functions $F$ from $H^2$ to
${\C}$ such that, for all multi-indices $(\alpha , \beta)$ in ${\cal M}_m$, the derivative
$\partial _q^{\alpha} \partial _p^{\beta }  F$ is well defined, continuous and bounded, and satisfies
%----
\be\label{3.1} | \partial _q^{\alpha} \partial _p^{\beta }  F (q , p)| \leq M \prod \varepsilon_{jk} ^{\alpha _{jk} + \beta _{jk} }.\ee
%-----
We set $S_m({\cal B} , \varepsilon)= \bigcup_{M\geq 0} S_m({\cal B} , M, \varepsilon)$.
For all $F\in S_m(\varepsilon)$, we set
$||F||_{m,\varepsilon}= \inf\{ M\geq 0 : F \in  S_m({\cal B}, M,\varepsilon)\}$.
If $F$ depends on one or several parameters, we say that $F$ is bounded in $S_m({\cal B} , \varepsilon)$ if the norm $||F||_{m,\varepsilon}$
is bounded independently of these parameters.

 We denote by $S_{\infty}({\cal B} , \varepsilon)$
the intersection of the classes $S_m({\cal B} , \varepsilon)$.
We say that $F$, depending on some parameters,
is bounded in $S_{\infty}({\cal B} , \varepsilon)$ if, for all $m$, the norm $||F||_{m,\varepsilon}$ is bounded independently of these parameters (but usually not on $m$).

We denote by $S_{\infty}^{mat}({\cal B} , \varepsilon)$ the analogous space for mappings taking values in
 ${\cal L} ( ({\C}^2)^{\otimes N})$.
\end{defi}

The Weyl calculus in infinite dimension  \cite{A-J-N} and \cite{A-L-N} enables to associate quadratic forms on the space ${\cal D}$ of Definition \ref{d2.2} with suitable functions. It is sufficient that the function satisfies   hypotheses $(H_1)$ and $(H_2)$ of \cite{A-L-N}
 (Section 1). The construction of the quadratic form is given in \cite{A-J-N} and, in a simpler way, in \cite{A-L-N} (Theorem 2.2). If the function $F$ belongs to $S_{2} ({\cal B} , \varepsilon)$,
 and if the sequence $\varepsilon_{mn}$ is summable then hypotheses $(H_1)$ and $(H_2)$ of \cite{A-L-N} are satisfied. In this case, one shows in Theorem 1.4 in \cite{A-J-N}, using Proposition 8.4 in \cite{A-J-N}, that  the  quadratic form initially defined  on ${\cal D}$  is the quadratic form of  a bounded operator in
 $L^2 ( B , \mu _{B , h/2})$ denoted $Op_h^{weyl} (F)$.

Besides the sets $S_{2} ({\cal B} , \varepsilon)$, there is another important class of functions satisfying  hypotheses $(H_1)$ and $(H_2)$ of \cite{A-L-N} and to which we may associate a quadratic form  on ${\cal D}$. These are  symbols  $L(q, p)$ which are continuous linear forms on $H^2$.
One shows (\cite{A-L-N}, Proposition 2.4) that the Weyl quadratic form   associated with such a symbol is the one of an operator
 from ${\cal D}$ to ${\cal D}$,  also denoted $Op_h^{weyl}(L) $. We may write $L(q , p) = (a , q) + (b , p)$,
with $a$ and $b$ in $H$. The Weyl quantized operator  $Op_h^{weyl}(L) $ associated
with this function is usually denoted $\Phi _S (a+ib)$ and called
Segal field.
Namely, using the  canonical Segal isomorphism
 ${\cal J}_h : {\cal F}_s (H_{\C}) \rightarrow L^2 (B, \mu _{B , h/2})$, this operator becomes,
 %----
\be\label{Segal_iso} {\cal J}_h ^{-1} Op_h^{weyl} (L)  {\cal J}_h = \sqrt {h}  \Phi_S (a+i b),\ee
 %----
 where the Segal field $\Phi_S (a+i b)$ is defined in  \cite{RE-SI-2}.
It may be defined without
 pseudodifferential formalism on  Fock space,
 using creation and annihilation operators. For any $a$ in $H$,
we denote by $Q_h (a)$ and $P_h (a)$ the operators associated with the functions $(q , p) \mapsto (a , q )$ and $(q , p) \mapsto (a , p )$ by the mapping
$Op_h^{weyl} $.

 For any $X=(a , b)$ in $H^2$ and for each $h>0$,  $\Psi_{X ,h}$ denotes the corresponding coherent state being an element of ${\cal F}_s (H_{\C})$ defined by,
%----
\be\label{csf} \Psi_{(a , b) , h}  = \sum _{n\geq 0}
 {e^{-{|a|^2+ |b|^2  \over 4h}} \over (2h)^{n/2} \sqrt {n!} } (a+ib) \otimes \cdots \otimes (a+ib).\ee
%----
In the space  $ L^2(B , \mu_{B , h/2})$, the coherent state is also denoted by
$\Psi_{(a , b) , h}$ and is defined by (\cite{A-L-N}),
%-----
\be\label{csl}\Psi_{a, b , h} (u) = e^{{1\over h} \ell _{ (a+ib)} (u)  -{1\over 2h}|a|^2 -
{i\over 2h} a\cdot b},\quad   {\rm a.e.}\  u\in B.\ee
%---
See Section \ref{s2.B} for  $\ell _{ (a+ib)}$.

\begin{prop}\label{p3.2} If $F$
belongs to $S_{\infty} (M, \varepsilon)$ and $G$ to $S_{\infty}  (M', \delta)$ then the
(possibly matricial) product $FG$ lies in $ S_{\infty}  (M M' ,  \varepsilon  + \delta )$.
 If $F$ takes scalar values, if $L $ is a continuous linear form on $H^2$ and if the sequence $(\varepsilon _{mn})$ is square
summable, then  the Poisson bracket $\{ F, L \}$ lies in $ S_{\infty}  (M'' ,  \varepsilon ) $
with,
%----
\be\label{3.2} M'' = \sum _{mn} \varepsilon _{mn} ( |L ( f_{mn}, 0)| + |L ( 0, f_{mn})| ).\ee
%----------
\end{prop}

With a given quadratic form  $Q$
on the space ${\cal D}$, one may define its Wick symbol by,
%----
\be\label{3.3}\sigma_h^{wick} (Q) (X) = Q ( \Psi_{X , h} , \Psi_{X , h} ),\qquad X \in H^2 \ee
%---
where the $\Psi_{X , h}$ are coherent states whose definition is  reminded in (16)
or in (17) of \cite{A-L-N}. One also defines the bi-symbol by,
%---
\be\label{3.4} (S_hQ) (X , Y) = { Q ( \Psi_{X , h} , \Psi_{Y , h} )  \over  <\Psi_{X , h} , \Psi_{Y , h}> },\qquad (X, Y) \in H^2 \times H^2.\ee
%----
 One similarly defines
$\sigma_h^{wick} (A)$ and $S_h A$ when $A$ is an operator from ${\cal D}$ to $L^2(B , \mu _{B , h/2})$.

\subsection{The functor $\Gamma $ and ${\rm d}\Gamma $. Photons number and energy.}\label{s3.B}

Let  $T$ be a selfadjoint operator, bounded or unbounded with domain $D(T)$ in $H^2$ and being
${\C}$-linear when $H^2$ and $H_{\C}$ are identified (that is to say, when
the operator $T$ commutes with the map ${\cal F}$ defined by ${\cal F} (q , p) = ( -p, q)$). One sets, for any finite sequence  $(u_1 , .. u_m)$
in $D(T)$,
%---
\be\label{3.5} d \Gamma (T) ( u_1 \otimes \cdots \otimes u_m) = (Tu_1)\otimes \cdots \otimes u_m + \cdots
 + u_1 \otimes \cdots \otimes (Tu_m ).\ee
 %-------
 By linearity, one defines an operator on the subspace of  ${\cal F}_s (H_{\C})$
generated by this type of elements. One shows (\cite{RE-SI-1}) that this operator is essentially selfadjoint and ${\rm d}\Gamma (T)$
   also  denotes its selfadjoint extension.

Let $S$ be a selfadjoint operator,  bounded in $H^2$, and commuting with ${\cal F}$.
Then the Wick symbol  of ${\rm d}\Gamma (S)$ is   $F(z) = < S(z), z> /2h$ if we identify
$(q , p)\in H^2$ with $q+ip\in H_{\C}$ (notation (\ref{1.4})).  This property enables to consider ${\rm d}\Gamma (S)$
as a quadratic Hamiltonian of the type studied in \cite{LA} by other methods.

To any operator $M$,
bounded in $H_{\C}$ with a norm smaller or equal than $1$, one associates an operator $\Gamma (M)$,
  bounded in ${\cal F}_s(H_{\C})$, such that, for all $u_1$, \dots,  $u_m$
   in $H_{\C}$, one has,
 %---
\be\label{3.6} \Gamma (M) ( u_1 \otimes ... \otimes u_m) = (Mu_1)\otimes ... \otimes (Mu_m).\ee
 %-------
In view of the usual isomorphism, these definitions are moved to the space $L^2(B , \mu_{B , h/2})$.

When $T= I$, the operator $N= {\rm d}\Gamma (I) $ is the number operator. When $T= h M_{\omega} $,
   where $M_{\omega} $ is the multiplication,
   in $H_{\C}$ or in $H^2$, by $\omega (k) = |k|$, the operator ${\rm d}\Gamma ( h M_{\omega})$, denoted $H_{ph}$, is
   the photons Hamiltonian.
If $M$ is the multiplication by $e^{it h \omega (k)}$,
then the  operator $\Gamma (M)$ is unitary in $H_{\C}$ and is precisely $e^{it H_{ph}}$.
For any rotation $R$ in $SO(3)$, one may define an operator $\pi (R)$ in $H$,
and then an operator   $\Pi (R)$ in  ${\cal H}_{ph}$ by,
%---
\be\label{3.7} (\pi (R) f ) (k) = R f ( R^{-1} k),\qquad \Pi (R)= \Gamma ( \pi (R)).\ee
%---
Thus, one defines a unitary representation  of $SO(3)$ in ${\cal H}_{ph}$ which is used in  Section \ref{s4}.

In the literature, the operator $Op_h^{weyl } (F_{a b} )$ is usually often denoted by
 $h^{1/2} \Phi_S (a+i b)$, where $F_{a b} $ is the linear form on $H^2$ defined
by $F_{a b} (q , p) = a \cdot q + b\cdot p$, with $(a,b)\in H^2$.

We denote the domain
of the operator $N^{m/2}$ by $W_m$ (Sobolev space) and it is equipped with,
%-----
\be\label{3.8} ||u||_{W_m}^2 = <(I +2h  N )^m u , u >.\ee
%---
In particular, for $m=1$, one has,
%-----
\be\label{3.9} ||u||_{W_1}^2= ||u||^2 + \sum_{(m,n)\in\Gamma}
||(Q_h(f_{mn}) + i P_h(f_{mn}))u||^2.\ee
%----

One knows that, when $L$ is a continuous linear form on $H^2$, the operator $Op_h^{weyl} (L)$
 initially defined from ${\cal D}$ into ${\cal D}$, is extended as an operator
from $W_1$ to $L^2(B , \mu _{B , h/2})$. See Lemma 2.3 in \cite{DG}.
 One also shows   in Proposition 2.8 of \cite{A-L-N} that, if $F$ belongs to $S_{\infty} ({\cal B}  ,  \varepsilon)$ with a summable sequence
  $(\varepsilon_{mn})$, then the  operator $Op_h^{weyl} (F)$ is also bounded in $W_m$.

The following  Proposition underlines the relation between the functor $\Gamma$ and  analogue of metaplectic operators, and between ${\rm d}\Gamma $ and quadratic Hamiltonians.
See Derezi\'nski-Gerard Lemmas 2.1 and 2.3 of \cite{DG} for an analogous statement.
See B. Lascar \cite{LA}  for an analogue of the metaplectic group in infinite dimension.
See Combescure Robert \cite{C-R} for quadratic Hamiltonians in finite dimension.

\begin{prop}\label{p3.3} Let $T$ be an unitary map in $H^2$  being
${\C}-$linear when identifying $H^2$ and $H_{\C}$ (in other words,
commuting with the map ${\cal F}$ defined by ${\cal F} (q , p) = ( -p, q)$). (These
hypotheses imply that $T$ is symplectic).
For any $X$ in $H^2$, one has,
%----
\be\label{3.10} \Gamma (T) \Psi_{X, h} = \Psi _{TX , h},\ee
%----
where the $\Psi_{X, h}$ are  coherent states. The Wick symbol of the operator
$\Gamma (T) $ is, with notation (\ref{1.4}) and identifying $X = (q , p) \in H^2$ with $ z = q+ip \in H_{\C}$,
%---
\be\label{3.11} <  \Gamma (T) \Psi_{z, h} , \Psi_{z, h}   > =  e^{-{1\over 2h } |z|^2 + {1\over 2h }
< Tz , z> }.\ee
%----
Let $\Phi $ be a continuous linear form on $H^2$. One has,
%----
\be\label{3.12}  (\Gamma (T)) ^{-1} Op_h^{weyl } (\Phi ) \Gamma (T) = Op_h^{weyl } (\Phi \circ T).\ee
%---
Let $S$ be a bounded  selfadjoint operator in $H^2$  commuting with ${\cal F}$, one has,
%----
\be\label{3.13}[ {\rm d}\Gamma (S) ,  Op_h^{weyl } (\Phi ) ] = i Op_h^{weyl } (\Phi \circ S \circ {\cal F}).\ee
%----
This commutator is a priori defined as a mapping from $W_3$ into $L^2(B , \mu _{B , h/2})$.\newline
If $S$ is unbounded then the result remains valid when the mapping $\Phi \circ S :
D(S) \rightarrow H^2$ is extended to a  bounded operator in $H^2$.
\end{prop}

The right hand side of  (\ref{3.13}) may be used to define a Poisson bracket.

We apply Proposition \ref{p3.3} with the map $T= \chi_t : H^2 \rightarrow H^2$ which is
identified to the multiplication by  $e^{-i t |k|}$ when  $H^2$ and $H_{\C}$ are identified.
That is to say, with $\omega (k) = |k|$,
%----
\be\label{3.14} \chi_t (q , p) = (q_t , p_t), \qquad
\left \{ \begin{array}{l}
q_t (k) = \cos (t \omega (k)) q(k) +  \sin (t \omega (k)) p(k)\cr  \cr
 p_t (k) = -  \sin (t \omega (k)) q(k) +  \cos (t \omega (k)) p(k)\cr
 \end{array} \right..
\ee
%----
In view of Proposition \ref{p3.3},
%---
\be\label{3.15}e^{-i{t\over h} H_{ph} }  \Psi_{q , p, h}  =  \Gamma (\chi_t ) \Psi_{q , p, h}  = \Psi_{\chi_t (q, p), h}.\ee
%----

We denote by   $H^{\omega }$  the set of  all $a$ in $H$ such that
$ |k| a(k)$ belongs to $H$ and by $H_{\omega }$  the set of all $a$ in $H$ such that
$  a(k)/ \sqrt {|k|}$ lies in $H$.

\begin{prop}\label{p3.4}  i) The space ${\cal D} \bigcap D(H_{ph})$ is dense in $W_1$ and in $D(H_{ph})$.

ii)  For all $a$ and $b$ in
$H$, let $F_{a b} $ be the linear form on $H^2$ defined
by $F_{a b} (q , p) = a \cdot q + b\cdot p$. Then, for all $a$ and $b$
in $H_{\omega }$, the operator $Op_h^{weyl } (F_{a , b}) $, initially defined
in ${\cal D}$,  is extended
as continuous  operator from $D( H_{ph})$ in $L^2(B , \mu _{B , h/2})$.
\end{prop}

 {\it Proof. i) } For all finite sequences $a_1,\dots,a_n$ in $H^{\omega }$,
and for all polynomial functions  $\Phi$, the function $f$ defined almost everywhere  on $B$ by,
%-
\be\label{3.16} f(x) = \Phi (\ell _{a_1} (x), \dots, \ell_{a_n} (x)),\quad
 {\rm \ a.\ e. \ } x\in B \ee
%---
 belongs to ${\cal D} \bigcap D( H_{ph})$ and the space of this kind of functions  is dense in $D( H_{ph})$.
Indeed, the set of functions $f$ defined in $(\ref{3.16})$ is, using Segal isomorphism, the space spanned by  symmetrized tensorial products
 $a_1\odot \cdots \odot a_n$, $a_j\in H^\omega$, which is a core
of $H_{ph}$ (\cite{RE-SI-2}).
Let us verify that the set of functions written as $(\ref{3.16})$ with the $a_k$ in $H$,
is dense in $W_1$. Indeed, one knows  that the norm of $W_1$ may also
be written as, using the number operator  $N$ or orthogonal projections  $f_n$
of an element $f$ of the Fock space into  $n$ photons spaces,
%---
\be\label{3.17}\Vert f \Vert _{W _1} ^2 = < (I + 2h N )f , f> = ||f||^2 +2h \sum _{n=1}^{\infty} n  ||f_n||^2.\ee
%----
Consequently, for any  $f$ in $W_1$, the sum $f_0 +\cdots + f_n$ tends to $f$ in $W_1$.
 Via the Segal isomorphism, this amounts to say that the set of functions under the form $(\ref{3.16})$, with  the $a_k$ in $H$ instead of  $H^{\omega}$, is dense in $W_1$. A function as in $(\ref{3.16})$ is a finite linear combination  of functions
under the form $g = : \ell _{a_1} \dots \ell _{a_n}: $ with the $a_k$ in $H$. One uses here the notation for the
Wick product of the functions $\ell _{a_k}$, that is to say, the range by the Segal  isomorphism of the symmetrized
 $a_1 \odot \cdots \odot a_n$ (\cite{J}).  Let $a_k^{(\nu)}$ be a sequence in $H^{\omega}$ converging to
 $a_k$ in $H$.  Let
$ g^{\nu}  = : \ell _{a_1^{\nu} } \dots \ell _{a_n^{\nu} } :$.
%----
The sequence of  functions $ g^{\nu}$ belongs to ${\cal D} \bigcap Dom (H_{ph})$ and the sequence $ (g^{\nu})$ tends to $g$ in $W_1$.  Using the Segal isomorphism, this amounts to say that the
sequence  $a_1^{(\nu)} \odot \cdots \odot a_n^{(\nu)}$ tends to
$a_1 \odot \cdots \odot a_n$ in $W_1$, which is standard.

ii) We use the following well-known inequalities, for all $\varphi\in D(H_{ph}^{1/2})$,  $a$ and
$b$ in $H_{\omega}$,
%----
 $$||Op_h^{weyl } (F_{a b} ) \varphi||\leq  C ( |a/\sqrt{\omega}| +|b/\sqrt{\omega}|  )
   ||H_{ph}^{1/2}\varphi|| + C  h^{1/2} (|a| + |b|)  ||\varphi||.
 $$
 %---
See \cite{B-F-S}.  Consequently,  for all $\varepsilon >0$, there exists $C_\varepsilon$  also depending on $a$, $b$ and $h$, such that,
 %---
 $$||Op_h^{weyl } (F_{a b} ) \varphi||\leq \varepsilon ||H_{ph}\varphi|| + C_\varepsilon  ||\varphi||,\qquad \varphi\in D(H_{ph}).  $$
 %--------

\hfill$\Box$

Note that $e^{-i{t\over h} H_{ph} }$ maps the space ${\cal D}$ of Definition \ref{d2.2} into itself.
Indeed, let $f$ belongs to ${\cal D}_E$
where $E$ is finite $n-$dimensional subspace of $H$. The function $X \mapsto < f, \Psi_{Xh}>$ lies in $L^1(E^2, \lambda)$ where $\lambda $ is the Lebesgue measure
and we may write the vector values integral,
%---
$$ f = (2\pi h)^{-n} \int _{E^2} < f, \Psi_{Xh}> \Psi_{Xh} d\lambda (X).$$
%----
Consequently, from (\ref{3.15}),
%---
$$ e^{-i{t\over h}H_{ph} } f = (2\pi h)^{-n} \int _{E^2} < f, \Psi_{Xh}>
\Psi_{\chi_t (X), h} d\lambda (X).$$
%----
Let $F_1$ and $F_2$ be the projections of the space $\chi_t(E^2)$ on the two components of
$H^2$. One sees that $e^{-i{t\over h}H_{ph} } f $
belongs to ${\cal D} _{F_1 + F_2}$.

Let $A$ be an  operator from  ${\cal D}$ into ${\cal D}$.   One defines another operator from ${\cal D}$ in ${\cal D}$ by,
\be\label{3.18} A^{free}(t)  =  e^{i{t \over h} H_{ph} } A e^{-i {t \over h} H_{ph} },\ee
for all
$t\in \R$.
%---
Then one has,
%---
\be\label{3.19} \sigma _h^{wick } ( A^{free}(t)) (q , p) =
\sigma _h^{wick } ( A) (\chi_t (q , p)).\ee
%---
  Indeed, for all $X$ in $ H^2$,
%--
$$ <  A^{free}(t) \Psi_{Xh} , \Psi_{Xh}  >  = < A  e^{-i  {t \over h} H_{ph} } \Psi_{Xh} ,
e^{-i {t \over h} H_{ph} }\Psi_{Xh} )= <A  \Psi_{\chi_t (X), h}, \Psi_{\chi_t (X), h} >
= \sigma _h^{wick } ( A) (\chi_t (X)).$$
%---

 Note that the symplectic maps analogous to $\chi _t$ and the corresponding unitary
 operators are used by B. Lascar in  Proposition 2.5 of \cite{LA}.

 We now consider Weyl symbol.

\begin{prop}\label{p3.5}   i) Let $F$  lying in $S_2 ({\cal B},  \delta)$ where the sequence
$(\delta _{mn})$ is summable. Set $t\in \R$. Suppose that $F\circ \chi_t$ lies in $S_2 ({\cal B},  \delta)$. Then, we have,
 %---
\be\label{3.20}  Op_h^{weyl} (F \circ \chi_t ) = e^{i{t \over h} H_{ph}}   Op_h^{weyl} (F)   e^{-i{t \over h}H_{ph}}.\ee
%---

 ii) Let $F$ and $G$ belong to $S_2 ({\cal B}, \delta)$ where the sequence
$(\delta _{mn})$ is summable. Fix $t\in \R$ and assume that,
%---
\be\label{3.21} Op_h^{weyl} (G)  =  e^{i{t \over h}H_{ph}} Op_h^{weyl} (F)   e^{-i{t \over h}H_{ph}}.\ee
%---
Then $G = F \circ \chi_t $.

iii) If $F$ is a continuous linear form on $H^2$ then we have (\ref{3.20}) for any $t\in \R$.
\end{prop}

{\it Proof.} i) Let $A_1$ and $A_2$ be the operators in the left and right hand sides of equality (\ref{3.20}). Denote by $\Phi_j = S_h (A_j)$ their Wick bi-symbol defined in (\ref{3.4}).  For any symbol $G$, also denote by $\Phi_G$ the Wick bi-symbol of $Op_h^{weyl}(G)$ similarly defined.
One then has $\Phi_1(X,Y)=\Phi_{F\circ \chi_t}(X,Y)$ and, from (\ref{3.10}) and since $\chi_t$ is unitary,  $\Phi_2(X,Y)=\Phi_{F}(\chi_t(X),\chi_t(Y))$, for $(X,Y)\in H^2$.
In view of the definition in Theorem 2.2 of \cite{A-L-N}, with functions $\ell $ defined
   in  Section \ref{s2.B} and the notion of stochastic extension in  Definition 4.4 of \cite{A-J-N} recalled
   in Section \ref{s2.B},
 %---
 $$ \Phi_F(X,Y)  = \int _{B^2} \widetilde  F(Z)
e^{{1 \over  h}( \ell _X ( \overline Z)  + \ell _{\overline Y}  ( Z ) - < X ,  Y >  } d\mu _{B^2, h/2} (Z),\qquad (X , Y) \in (H^2 \times H^2)$$
%------
where $\widetilde  F$ denotes the stochastic extension of $F$ which exists from  Proposition 8.4 of
\cite{A-J-N}.   In the above exponent, we have identified $H^2$ with $H_{\C}$ and the notation $< \ , \ >$ is the one of (\ref{1.4}). Let us check  that $\Phi_1=\Phi_2$   on $ H^2$. To this end,  set for all $X$, $Y$
and $Z$ in $H^2$,
%----
$$ G_2  (Z) = F(Z) e^{{1 \over  h}( < \chi _ t (X ) ,   Z>
+ < Z ,    \chi _t (Y ) >  ) - < \chi _t (X) ,  \chi_t(Y) >   }$$
%----
and $G_1 (Z) = G _2( \chi_t (Z))$. These two functions of $Z$ have stochastic extensions
 $\widetilde G_2 $ and $\widetilde G_1 = \widetilde {G_2 \circ \chi_t }$ in $L^1(B^2,\mu_{B^2,h/2})$.
One has, for all $X$ and  $Y$   in $H^2$, since $\chi_t$ is unitary
%----
$$ \Phi _2 (X , Y) =  \int _{B^2}  \widetilde G_2  (Z) d\mu _{B^2, h/2} (Z),\qquad \Phi _1  (X , Y) =  \int _{B^2} \widetilde {G_2 \circ \chi_t }
 (Z) d\mu _{B^2, h/2} (Z).  $$
 %---
These integrals are equal from Proposition \ref{p2.3}. Consequently, one indeed has $\Phi_1=\Phi_2$   on $ H^2$. Therefore, these two operators $A_1$ and $A_2$ have the same Wick bi-symbol. One concludes according to Lemma 2.7
 of \cite{A-L-N}.

ii) Since $F$ and $G$ lies in $S_2 ({\cal B}, \delta)$ then they have stochastic extensions $\widetilde F$
and $\widetilde G$ in $L^2 (B , \mu _{B , h/2})$. From Proposition \ref{p2.3}, this holds true for $F \circ \chi_t$, and
$F \circ \chi_t$ satisfies hypotheses $(H_1)$ and $(H_2)$ in Theorem 2.2 of \cite{A-L-N}. Then, one may associate a quadratic form  $Q_h^{weyl} (F \circ \chi_t)$, and a priori not an operator, with
$F \circ \chi_t$,  and similarly for
 $G$. The same method as in  i)  shows that $Q_h^{weyl} (G-F\circ \chi_t)(f,g)=0$ for $(f,g)\in {\cal D}^2$.
 This implies that $(G-F\circ \chi_t)\circ P_E$ is vanishing when choosing
$(f,g)\in {\cal D}_E^2$, for all vectors subspaces $E$ in $H$, and therefore that $G-F\circ \chi_t = 0$.

iii) If $F$ is a continuous linear form on $H^2$ then this also true for $F \circ \chi _t$.
The functions $F$ and  $F \circ \chi _t$ remain fixed when applying the operator $H_{h/2}$ defined in (37) of \cite{A-L-N}. Consequently, from the equality (38) of \cite{A-L-N},
 Weyl and  Wick symbols of these operators are equal.  Point iii) then follows from (\ref{3.13}).

\hfill \carre

\subsection{Field operators.}\label{s3.C}

For $x$ in $\R^3$, we now recall  the standard definitions  of the  six operators $B_j (x)$ and
$E_j(x)$ ($1 \leq j \leq 3)$, unbounded in the  space ${\cal H}_{ph}$,
corresponding to the three components of the magnetic field and to the three components of the electric field. These operators are usually defined with the Segal fields
formalism (see, e.g., \cite{F-1}\cite{F-2}\cite{B-F-S}\cite{SP}). However, to use
pseudodifferential calculus, we associate these  operators  with symbols noted $B_j(x, q,
p)$ and $E_j(x, q, p)$ which are linear in $(q, p)$.
To define  symbols $B_j(x, q, p)$,  one chooses a function $\chi$ in ${\cal S}(\R)$ vanishing in a neighborhood  of $0$ and
set,
%----
\be\label{3.22}B_{j,x}(k) = {i\chi(|k|)|k|^{1\over 2} \over (2\pi)^{3\over 2}}
e^{-ik\cdot x } {k\wedge e_j \over |k|},\qquad k\in\R^3\backslash\{0\} \ee
%----
where $(e_1,e_2, e_3)$ is the canonical basis  of $\R^3$.

The function $\chi$ is assumed to vanish in a neighborhood of the origin in order to apply the results in \cite{A-J-N} and \cite{A-L-N} in their current forms. A possible improvement of these results could allow to avoid this hypothesis.

Next, we set,
%----
\be\label{3.23}B_j (x, q, p) = ({\rm Re} B_{j,x} ,  q ) +  ( {\rm Im }   B_{j,x},  p ).\ee
%------

Explicitly,  one has for instance,
%---
$$ B_3^{free} (x , t, q, p) = (2\pi)^{-{3\over 2} }\int _{\R^3} \chi(|k|)|k|^{1\over 2}
    \Big [ - \sin   ( k\cdot x - t|k|) (k_1 q_2(k) - k_2 q_1(k)) $$
    %---
$$- \cos   ( k\cdot x - t|k|) (k_1 p_2(k) - k_2 p_1(k))\Big ] { dk \over |k|}. $$
%--

We denote by $B_j (x)$ the unbounded operator whose symbol is $ B_j (x, q, p)$.

We  also use the helicity operator  $J$ from $H^2$ to $H^2$ and defined by,
%----
\be\label{3.24}J(q, p) (k) = \left ( {k\wedge q(k) \over |k| } , {k\wedge p(k) \over |k| } \right ),\qquad
k\in \R^3 \setminus \{ 0 \}.\ee
%-----
One then sets,
%----
\be\label{3.25}  E_j (x, q, p) =  - B_j (x , J(q, p)).\ee
%---

We use the  commutation relations below, where we set,
%---
\be\label{3.26}\rho (x) = (2 \pi )^{-3} \int _{\R^3 } |\chi (k)| ^2 \cos (k \cdot x ) dk.\ee
%------
One has,
%----
\be\label{3.27} [   E_1(x), B_2(y) ] = i h  (\partial_3 \rho ) (y - x).\ee
%------
 The other commutators follow by anti-symmetry and circular permutations. Also,  $ [ E_j(x), B_j(y) ]  = 0$.
Moreover,
%---
\be\label{3.28}[B_j (x), B_m (y) ] =  [E_j (x), E_m (y) ] = 0.\ee
%-----------

Indeed, if $A$ and  $B$ are two operators with linear symbols $F$
 and $G$ (two Segal fields) then one knows that the commutator $[A , B]$ is equal to the constant
$ (h/i ) \{ F , G\}$ where $\{ F , G\}$ is the Poisson bracket. When $F(q, p)
= (a , q) + (b , p)$ and $G(q, p) = (a', q) + (b', p)$, one has $\{ F , G\} = (b, a') - (a, b')$.

Next, the fields free evolution operators are defined by,
\be\label{3.29} B_j^{free}(x, t) = e^{i{t\over h}  H_{ph}} B_j(x) e^{-i{t\over h}  H_{ph}},\quad
   E_j^{free}(x, t) = e^{i{t\over h} H_{ph}} E_j(x) e^{-i{t\over h} H_{ph}}.\ee
 %-----
From Proposition 3.4 (point iii), these operators have Weyl and Wick symbols,
 %----
\be\label{3.30}B_j^{free}(x,t,q,p) = B_j(x,\chi_t(q,p)),\qquad
 E_j^{free}(x,t,q,p) = E_j(x,\chi_t(q,p)).\ee
 %----

The operators valued vector fields
${\bf B}^{free}(x, t) = (B_1^{free}(x, t),B_2^{free}(x, t),B_3^{free}(x, t))$
and ${\bf E}^{free}(x, t)$ satisfy
Maxwell equations in vacuum,
\be\label{3.31}{\rm div}\, {\bf B}^{free}(x, t) = {\rm div}\, {\bf E}^{free}(x, t) = 0\ee
\be\label{3.32} {d \over dt} {\bf B}^{free}(x, t) = -{\rm rot}\, {\bf E}^{free}(x, t),\qquad
   {d \over dt} {\bf E}^{free}(x, t) =  {\rm rot}\, {\bf B}^{free}(x, t).\ee
The symbols of these operators satisfies the same equations.

\section{Spin-photon interaction model.}\label{s4}

The below model  describing the quantized electromagnetic field
interacting with $N$ fixed  $1/2$ spin particles may be found in
Reuse \cite{REU} (Section 4.11) in the case of one  particle and then extended in a straightforward way for $N$ particles.

The Hilbert space of the quantized field  is,
%---
\be\label{4.1} {\cal H}_{ph} = L^2 ( B , \mu _{B , h/2}) \simeq {\cal F}_s (H_{\C})\ee
%---
where $H$ is the Hilbert space of Section \ref{s2.A}, $B$ and $\mu _{B , h/2}$  are the Banach space
and the measure defined in Section \ref{s2.B}. In this space ${\cal H}_{ph}$, one considers
the photons Hamiltonian,
%---
\be\label{4.2} H_{ph} = h {\rm d}\Gamma (M_{\omega})\ee
%---
($M_\omega$ being the multiplication by $|k|$).

The Hilbert space of a $1/2$
spin  particle, fixed at the origin, is $\C^2$. The
Hamiltonian of this particle in an exterior magnetic field
${\bf \beta } = (\beta_1,\beta_2,\beta_3) \in \R^3$ is,
\be\label{4.3} H_{mag } = \sum_{j=1}^3 \beta_j \sigma_j\ee
 where the
matrices $\sigma_j$ ($1\leq j \leq 3$) are the Pauli matrices,
%----
\be\label{4.4}
\sigma_1 = \begin{pmatrix}  0 & 1 \\ 1 & 0   \end{pmatrix}
\qquad
   \sigma_2 = \begin{pmatrix} 0 & -i \\ i & 0   \end{pmatrix}
 \qquad
   \sigma_3 = \begin{pmatrix}  1 & 0 \\ 0 & -1  \end{pmatrix}.\ee
%----

The Hilbert space of the entire system  constituted with the photons field and with $N$  fixed $1/2$ spin particles at the points in $\R^3$ denoted  $a_1,\dots,a_N$ is
${\cal F}_s(H)\otimes ( \C^2  )^{\otimes N} $.
The Hamiltonian of this whole system  is,
%----
\be\label{4.5} H (h) = H_0 + h H_{int},\qquad  H_0 = H_{ph}\otimes I\ee
%----
with,
%----
\be\label{4.6} H_{int} = \sum _{\lambda = 1}^N \sum_{j=1}^3
( \beta _j + B_j(a_\lambda))\otimes  \sigma _j ^{[\lambda]}\ee
%---
where, for all operators $A$ belonging to ${\cal L} ({\C}^2)$, $A^{[\lambda]}$ denotes
$I \otimes \cdots A \cdots \otimes I$ with $A$ located at the $\lambda ^{th}$ position.

\begin{prop}\label{p4.1} The operator $H(h)$ above, of domain $D(H(h)) = D(H_0) = D(H_{ph}) \otimes ( {\C}^2 )^{\otimes N}$,
is selfadjoint. For all $x$ in $\R^3$, the operators $B_j(x)$ and $E_j (x)$
($1 \leq j \leq 3$) are bounded from $D(H(h))$ into ${\cal H} _{ph} \otimes ( {\C}^2
 )^{\otimes N}$.
\end{prop}

Indeed, from  Proposition \ref{p3.4} point ii), one has for all $a$ and $b$
in the space $H_{\omega}$  defined before this Proposition,   $D(H_{ph})\subset D(Op_h^{weyl} (F_{a b}) )$ where
$F_{a b} (q , p) = a \cdot q + b \cdot p$. Since the operators $B_j(0)$ are of this kind, one sees,
for all $\phi\in D(H_{0})$ and for any $\varepsilon >0$,
 %---
\be\label{4.7}||H_{int} \phi||\leq \varepsilon ||H_{0}\phi|| + C_\varepsilon ||\phi||\ee
 %--------
and one concludes applying Kato-Rellich Theorem.

We can also consider the case of one $1/2$ spin  particle, fixed at the origin, interacting both with a constant magnetic field and with a non quantized magnetic rotating field in an orthogonal plane to the constant field, and also interacting with the quantized field. Denoting the constant field $\beta = (0, 0 , \beta _3)$ and setting
 $(B _1 \cos ( \omega  t) , B_1 \sin ( \omega t) , 0)$ as the rotating field
 then one gets the following Hamiltonian $H(h, t)$,
%----
\be\label{4.8}H(h, t) = H(h) +  h H_{spin} (t),\qquad  H_{spin} (t)  = ( B _1 \cos ( \omega  t) (I \otimes  \sigma _1) +  B_1 \sin ( \omega t)
( I \otimes \sigma _2)),\ee
%---
where $H(h)$ is the Hamiltonian in  (\ref{4.5}) with $\beta = (0, 0, \beta _3)$. The next Theorem verifies
that, using a suitable transform, one may reduce this model to a time independent model with an Hamiltonian similar to the one in (\ref{4.5}).

\begin{theo}\label{t4.2} One has,
%---
\be\label{4.9} P(t) ^{\star} \left ( ih  {\partial \over \partial t} - H(h, t)     \right )  P(t) =
ih  {\partial \over \partial t} - H_{TR},\ee
%---
where
%----
\be\label{4.10} H_{TR} = H_0 + \sum_{j=1}^3  B_j(0)\otimes  \sigma _j  + B _1 (I \otimes  \sigma _1)
+ ( \beta _3  - {\omega \over 2 }  ) (I \otimes  \sigma _3)\ee
%--------
and $ P(t) =  P_{ph} (t) \otimes  P_{spin} (t)$ with $P_{ph} (t)$ and $ P_{spin} (t)$
being operators in ${\cal H}_{ph}$ and ${\C}^2$ defined as following. One has,
$ P_{ph} (t) = \Pi (R (-\omega t)))$ where $\Pi $ is the $SO(3)$ representation
in ${\cal H}_{ph}$ defined in (\ref{3.7}) and
 $R( \omega t)$ is the rotation with angle $\omega t$ in the horizontal plane. Moreover,
%---
\be\label{4.12} P_{spin} (t) = \begin{pmatrix} \alpha (t) & 0\cr 0 & \alpha (-t) \end{pmatrix}, \qquad
\alpha (t) = e^{ -i {\omega t \over 2}}.\ee
%-----
\end{theo}

{\it Proof.} Clearly, $ P_{spin} (t) ^{-1} \sigma_j P_{spin} (t) = \sigma _j (t)$ with,
%----
\be\label{4.13} \sigma_1 (t) = \begin{pmatrix}  0 & e^{i\omega t} \cr  e^{-i\omega t}& 0 \cr \end{pmatrix},
\qquad \sigma_2 (t) = \begin{pmatrix} 0 & -ie^{i\omega t} \cr  i e^{-i\omega t} & 0 \cr \end{pmatrix},
\qquad\sigma_3 (t) = \sigma_3.\ee
%----
Consequently,
%---
$$ P_{spin} (t) ^{-1} \Big ( \cos (\omega t) \sigma_1 + \sin (\omega t) \sigma_2 \Big )
P_{spin} (t)  = \sigma_1,$$
%------
$$  P_{spin} (t) ^{-1} \sigma _3  P_{spin} (t) = \sigma_3, $$
%----
$$  P_{spin} (t) ^{-1} ih  {\partial \over \partial t}  P_{spin} (t) =
 ih  {\partial \over \partial t} + {h \over 2 } \omega \sigma_3. $$
 %----
Thus,
 %----
$$ P_{spin} (t) ^{-1} \left (   ih  {\partial \over \partial t}  - h
 \Big ( B_1 \cos (\omega t) \sigma_1 +B_1  \sin (\omega t) \sigma_2  + \beta _3 \sigma_3
 \Big ) \right ) P_{spin} (t)$$
 %---
\be\label{4.14} =  ih  {\partial \over \partial t}  -  h   \Big ( B_1 \sigma_1 + ( \beta _3
 - {\omega \over 2})  \sigma _3 \Big ).\ee
 %-----
One knows that $H_{ph} = d \Gamma ( h M)$ where $M$ is the multiplication by $|k|$ which
 commute with $U$ in such a way that $H_{ph} $ commutes with $\Gamma (U)$. Then,
 %---
\be\label{4.15} P_{ph} (t) ^{-1} H_{ph}P_{ph} (t) =  H_{ph}.\ee
 %---
 One now checks that,
 %----
\be\label{4.16} \sum _{j=1}^3  P (t) ^{-1}  \Big ( B_j (0) \otimes \sigma _j \Big ) P (t) = \sum _{j=1}^3
  B_j (0) \otimes \sigma _j.\ee
 %------
From Proposition \ref{p3.3}, one has,
 %---
 $$ \sum _{j=1}^3  P (t) ^{-1}  \Big ( B_j (0) \otimes \sigma _j \Big ) P (t) =
 \sum _{j=1}^3  Op_h^{weyl} ( C_j (t, \cdot )) \otimes  \sigma_j (t) $$
 %---
 where $\sigma_j (t) $ is defined in (\ref{4.13}) and
 %---
 $$ C_j (t, q, p)   = B_j (0, U(t) (q , p)). $$
 %---
A direct computation shows that,
 %---
 $$ C_1 (t, \cdot) + i C_2 (t, \cdot) = e^{i\omega t} ( B_1 (t, \cdot) + i B_2 (t, \cdot)),\qquad C_1 (t, \cdot) - i C_2 (t, \cdot)= e^{-i\omega t} ( B_1 (t, \cdot) - i B_2 (t, \cdot))
 $$
%---
and that $C_3 (t, \cdot) = B_3 (t, \cdot)$.  Equality (\ref{4.16}) is then obtained. Therefore, equality (\ref{4.9}) follows
from (\ref{4.14})(\ref{4.15}) and (\ref{4.16}).

\hfill$\Box$

\section{Proof of Theorem \ref{t1.1} $(i)(ii)(iv)$.}\label{s5}

We first begin by defining a numerical sequence  $(\varepsilon _{mn} (t))$ which, associated
with the basis $(f_{mn})$ in Section \ref{s2.C}, plays a special role in our class of  operators.

We set,
%----
\be\label{5.1}H^{free}_{int}(t) =  e^{i{t \over h} H_0} H_{int} e^{-i{t \over h} H_0}
= \sum _{\lambda = 1}^N \sum_{j=1}^3
( \beta _j + B_j^{free} (a_{\lambda} , t)\otimes  \sigma _j ^{[\lambda]}.\ee
%---
We use the basis $(f_{mn})$ defined in Section \ref{s2.C}. The operators $P_h(f_{mn})$ and $Q_h(f_{mn})$ denote
momentum  and position operators  corresponding to the element $f_{mn}$.

\begin{prop}\label{p5.1} Let $(f_{mn})$ be the basis chosen in Section \ref{s2.C}
and let $\chi$  be in ${\cal S}(\R)$ vanishing in a neighborhood of the origin. Then there exists a  sequence $\varepsilon_{mn}(t)$
satisfying,
%---
\be\label{5.2} \Vert [ P_h(f_{mn}) , H^{free}_{int}(t) ] \Vert + \Vert [ Q_h(f_{mn}) , H^{free}_{int}(t) ] \Vert
\leq h \varepsilon _{mn} (t).\ee
%----
Moreover, the sequence $\varepsilon_{mn}(t)$  is  rapidly decreasing and in particular, it is summable, and
$\varepsilon _{mn} (t)$ is a non decreasing function  of $|t|$.
\end{prop}

{\it Proof.}
Set,
%--------
\be\label{5.3} B_{jxt}(k) = {i\chi(|k|)|k|^{1\over 2} \over (2\pi)^{3\over 2}}
e^{i( t|k | - k\cdot x   )} {k\wedge e_j \over |k|},\quad k\in\R^3\backslash\{0\}.\ee
The element $B_{j,a_{\lambda } ,s}$ of $H$ has its three  components belonging to ${\cal S} (\R^3)$ and theirs norms
 $N_{\alpha \beta }$ of the Proposition \ref{p2.4} are bounded by  a constant $C_{\alpha \beta }(t)$ which  is bounded on every
 compact set in $\R$. From  Proposition \ref{p2.4}, one sees,
 %---
\be\label{5.4} m^{\alpha } n^{\beta } | (B_{j,a_{\lambda } ,s} , f_{mn})| \leq C_{\alpha \beta } (t).\ee
 %-
Thus, we may choose  a function $\delta_{mn}(t)$ being a non decreasing function of $|t|$, which is summable for all $t$
and satisfying,
 %------
\be\label{5.5}  \sup_{|s| < |t|}   \sup _{\lambda \leq N}  \left (  \sup_{j\leq 3}|(B_{j,a_{\lambda } ,s} , f_{mn})| \right )
\leq \delta_{mn}(t).\ee
 %----
In view of (\ref{3.30})(\ref{3.14})(\ref{3.22}) and (\ref{3.23}), $B_j^{free} (x , t, q, p)$ is the scalar product  of $B_{jxt}$
with $(q , p)$.  According to the expression  $H^{free}_{int}(t)$ in (\ref{5.1}) together with (\ref{5.5}), there exists
a constant $K>0$ such that,
%---
$$  \Vert [ P_h(f_{mn}) , H^{free}_{int}(t) ] \Vert + \Vert [ Q_h(f_{mn}) , H^{free}_{int}(t) ] \Vert
\leq K h \delta _{mn} (t).$$
%---
One therefore deduces the Proposition with $\varepsilon _{mn } (t) = K \delta _{mn } (t)$. Note that this family is summable
from (\ref{5.4}).

\hfill $\Box$

We  now consider the following
 operator
 (interaction picture) defined by (\ref{1.1}),
%----
\be\label{5.6}  U^{red}_h(t) =  e^{i{t \over h}  H_0} e^{-i {t \over h}  H(h)}.\ee
%----
In order to prove the pseudodifferential nature of $U^{red}_h(t)$, we now study iterated commutators of that operator with position $Q_h (f_{mn})$ and momentum $P_h (f_{mn})$  operators associated with the elements $f_{mn}$ of the basis ${\cal B}$ defined in Section \ref{s2.C}. We prove that, these iterated commutators, a priori defined as quadratic forms
 on some subspaces, are extended as bounded operators
in  ${\cal H}_{ph}\otimes ( {\C}^2 )^{\otimes N}$. Let us start with simple commutators.  For each continuous linear form $F$ on $H^2$,
one defines the following commutator as  a quadratic form  on $W_1\otimes ( {\C}^2 )^{\otimes N}$ by,
%----
$$[ Op_h^{weyl} (F )\otimes I ,  U^{red}_h(t) ] (f , g) = <  U^{red}_h(t) f, (Op_h^{weyl} (F )^{\star}\otimes I ) g >
 $$
%---
\be\label{5.7}  - < ( Op_h^{weyl} (F )\otimes I)  f, U^{red}_h(t) ^{\star} g>.\ee
%--------
If $A_1, \dots , A_m$ are operators, each being either one of the
$P_h ( f_{qr})\otimes I$ or one of the $Q_h ( f_{qr})\otimes I$ associated with one  of the element of the basis ${\cal B}$ (See sections \ref{s2.C} for this basis and \ref{s3.A} for position and momentum operators),   one similarly defines the iterated commutator, with the notation ${\rm ad} (A) B=[A,B]$,
%---
\be\label{5.8} {\rm ad} (A_1 )\dots {\rm ad} (A_m  )
U^{red}_h(t)\ee
%--
being a priori defined as a bilinear form  on $W_m\otimes ( {\C}^2 )^{\otimes N}$.
Let us now prove that these commutators are extended as bounded operators in ${\cal H}_{ph }
\otimes ( {\C}^2 )^{\otimes N}$. We begin with order one  commutators.

\begin{prop}\label{p5.2}   For all $f$ and $g$ in $W_1\otimes  ( {\C}  ^2 )^{\otimes N} $, for
any continuous linear form $F$  on $H^2$ and for all $t\in \R$, one has,
%---
\be\label{5.9}  [ Op_h^{weyl} (F )  \otimes I  ,  U^{red}_h(t) ] (f , g) = < C(t) f , g>,\ee
with,
%---
\be\label{5.10}  C(t) f = - i \int_0^t  U^{red}_h(t)  U^{red}_h(s)^{\star} [  Op_h^{weyl} (F) \otimes I ,
H^{free}_{int}(s) ] U^{red}_h(s) ds.\ee
%----
\end{prop}

{\it Proof. First step.}
We denote by $H_{\omega}$ the set of all $q$ in  $H$ such that
the function $k \mapsto q (k)/ |k|^{1/2}$ belongs to $H$. For any
$a$ and $b$ in $H$, set $F_{ab} (q , p) = (a , q) + (b , p)$.
We first prove equalities (\ref{5.9}) and (\ref{5.10}) when $f$ and $g$
lies in ${\cal D}\otimes ( {\C}^2  )^{\otimes N}  \bigcap D (H(h))$ and when $F = F_{ab}$ with $a$ and $b$
 in $H_{\omega}$. Set,
%---
$$ Z(t) = U^{red}_h(t)  ^{\star} (Op_h^{weyl} (F_{a  , b} ) \otimes I)   U^{red}_h(t)  f.$$
%---
This indeed defines an element  of ${\cal H} _{ph} \otimes ( {\C}  ^2  )^{\otimes N} $. Since
$ U^{red}_h(t) $ maps $D (H(h) )$ into itself from (\ref{5.6}) and  Proposition \ref{p4.1}, and since,
from  Proposition \ref{p3.4}, for all $a$ and $b$ in
$H_{\omega}$, the operator $Op_h^{weyl} (F_{a  , b} ) \otimes I$ maps $D (H(h) )$ into
${\cal H} _{ph}\otimes ( {\C}^2 )^{\otimes N}$. From (\ref{3.22}), $B_{j x}$ belongs to $H_{\omega}$,
 and  $H_{int}$ maps  $D (H (h) )$ into
${\cal H} _{ph} \otimes ( {\C}^2 )^{\otimes N}$ and we have, for all $f$ in $D (H (h) )$,
%---
\be\label{5.11} {d \over dt} U^{red}_h(t)  f = - i H^{free}_{int}(t) U^{red}_h(t)f.\ee
Consequently, if $g$ is also in $D (H (h) )$,
%---
$$ {d \over dt} < Z(t), g> = -i <H^{free}_{int}(t)  U^{red}_h(t)  f  , (Op_h^{weyl} (F_{a  , b} ) ^ {\star } \otimes I)  U^{red}_h(t)g > $$
%---
 $$ + i <  (Op_h^{weyl} (F_{a  , b} )\otimes I) U^{red}_h(t)  f  , H^{free}_{int}(t)  U^{red}_h(t)  g>      $$
%-----
$$ = - i < [  (Op_h^{weyl} (F_{a  , b} )\otimes I)  ,  H^{free}_{int}(t) ] U^{red}_h(t)  f  , U^{red}_h(t)  g >. $$
%------
The commutator $ [  Op_h^{weyl} (F_{a  , b} )\otimes I ,  H^{free}_{int}(t) ] $
is classically a constant matrix and then, is a bounded operator. We may then write,
%-----
$$ Z'(t) = - i U^{red}_h(t)^{\star}[  Op_h^{weyl} (F_{a  , b} )\otimes I ,  H^{free}_{int}(t) ] U^{red}_h(t)  f. $$
%----
As a consequence, for all $f$  in $({\cal D}\otimes  ( {\C}^2  )^{\otimes N} ) \bigcap D (H (h) )$ and for and
 $a$ and $b$  in $H_{\omega}$, we have,
%---
$$  U^{red}_h(t)  ^{\star} Op_h^{weyl} (F_{a  , b} )   U^{red}_h(t)  f = Op_h^{weyl} (F_{a  , b} )  f -
i \int _0^t
U^{red}_h(s)^{\star}[  Op_h^{weyl} (F_{a  , b} ) ,  H^{free}_{int}(s) ] U^{red}_h(s)  f ds
  $$
%---
 which implies, when acting  $  U^{red}_h(t)$ on the left side,
%---
\be\label{5.12} [  Op_h^{weyl} (F_{a  , b} ), U^{red}_h(t)  ] f = - i \int_0^t  U^{red}_h(t)
U^{red}_h(s)^{\star}[  Op_h^{weyl} (F_{a  , b} ) ,  H^{free}_{int}(s) ] U^{red}_h(s)  f ds.\ee
%---
The Proposition is then proved in this case.

{\it Second step.}  For any $(a , b)$ in $H^2$, we choose a sequence $(a_j , b_j)$ in
$H _{\omega }^2$ converging to $(a , b)$ in $H^2$. From a standard result recalled in
 Proposition 2.8 of \cite{A-L-N}, for all $f$ in
$ W_1 \otimes ( {\C}^2 )^{\otimes N} $,
%----
$$ \Vert  Op_h^{weyl} (F_{a  , b} - F_{ a_j , b_j}) f \Vert \leq C ( |a - a_j|+ |b- b_j|)
\Vert f \Vert _{W_1}.  $$
%---
Consequently,
%----
$$ \lim _{j\rightarrow \infty } \Vert  Op_h^{weyl} (F_{a  , b} - F_{ a_j , b_j}) f \Vert = 0,\qquad
\lim _{j\rightarrow \infty } \Vert [ ( F_{a  , b} - F_{ a_j , b_j}) \otimes I , H_{int}^{free } (t , \cdot ) ]
\Vert = 0$$
%-------
Therefore, the Proposition holds true for all $f$ and $g$ in
$({\cal D}\otimes   ( { \bf C} ^2  )^{\otimes N}  )\bigcap D (H (h) ) \subset  W_1\otimes ( {\C}  ^2 )^{\otimes N}  $, and for any
$a$ and $b$ in $H$.

{\it Third step.} Taking into account Proposition \ref{p3.4}, the space $({\cal D}\otimes  ( {\C} ^2 )^{\otimes N}  )\bigcap D (H (h) )$
is dense in $W_1 \otimes  ( {\C} ^2  )^{\otimes N} $ and the Proposition then follows.

\hfill \carre

We now turn to iterated  commutators. For all $(m, n)$, we set
 $\widetilde P_h ( f_{mn} ) = P_h ( f_{mn} ) \otimes I$ and for all multi-index
$(\alpha , \beta)$ we define a quadratic form  on
${\cal D} \otimes  ( {\C} ^2  )^{\otimes N} $ denoted   $({\rm ad} \widetilde P_h )^{\alpha} ({\rm ad} \widetilde Q_h )^{\beta}  U_{h}^{red}  (t)$,
 as in Section 1 of \cite{A-L-N}.

\begin{prop}\label{p5.3}  Let $U_h^{red}(t)$ be the operator defined in (\ref{1.1}) or (\ref{5.6}).
Then,  for all multi-index $(\alpha , \beta)$, the quadratic form
$({\rm ad} \widetilde P_h )^{\alpha} ({\rm ad} \widetilde Q_h )^{\beta}  U_{h}^{red}  (t)$
is the quadratic form of a bounded operator, written with the same notation, which satisfies,
%----
\be\label{5.13} \Vert  ({\rm ad} \widetilde P_h )^{\alpha} ({\rm ad} \widetilde Q_h )^{\beta}  U_{h}^{red}  (t)
\Vert \leq  \prod _{mn }
( h |t | \varepsilon _{mn} (t) )^{\alpha _{mn} + \beta _{mn} }\ee
%--
where $\varepsilon _{mn} (t) $ is defined in  Proposition \ref{p5.1}.
\end{prop}

{\it  Proof.}
Let $  A_1,\dots,   A_n$  be a finite sequence of operators, each operator $A_j$
being one of the
$P_h ( f_{qr})\otimes I$ or one of the $Q_h ( f_{qr})\otimes I$ associated with one of the element of the basis ${\cal B}$,
the corresponding indice $(q , r)$ being denoted $\psi (j)$.
Iterating Propositon \ref{p5.2}, one gets,
 %----
$$     {\rm ad} A_1 ...  {\rm ad} A_n U_h^{red} (t) =
 (-i)^n  \sum _{\varphi \in S_n}  \int _{  \Delta _ n(t)}
  U(t, s_n)[ A_{\varphi (n)} , H^{free}_{int}(s_n)]
  U(s_n, s_{n-1} )[ A_{\varphi (n-1)}, H^{free}_{int}(s_{n-1})]$$
  %---
\be\label{5.14} U(s_2, s_1)[ A_{\varphi (1)}, H^{free}_{int}(s_1)]
U(s_1, 0) ds_1 ...ds_n\ee
%---
 where we set $U(t , s) = U_h^{red} (t) (U_h^{red} (s) )^{\star}$, where
  $S_n$ is the set of bijections $\varphi $  in
the set $\{ 1, ... , n \}$, and, if $t>0$,
%----
\be\label{5.15} \Delta _ n(t) = \{ (s_1, ..., s_n) \in
\R^n , 0 < s_1 < ... < s_n < t  \} .\ee
%----
From  Proposition \ref{p5.1} and since the $\varepsilon _{mn} (t)$ is non decreasing, we obtain that, if $|s| < |t|$,
%---
$$ \Vert [ A_j  ,H^{free}_{int}(s)] \Vert  \leq  h  \varepsilon _{\psi (j)}(t). $$
%----
 Consequently, since $U(t, s)$ is unitary,
%----
$$ \Vert {\rm ad} A_1 ...  {\rm ad} A_n  U_h^{red} (t) \Vert
\leq  \sum _{\varphi \in S_n} \prod _{j=1}^n
h    \varepsilon _{ \psi ( \varphi (j))} (t)  \int_{\Delta_n (t)} ds_1\dots ds_n
 $$
%---
 where $\Delta _{n}(t)$ is defined in (\ref{5.15}).
The last integral equals to $|t|^n / n!$.
The factors in front  this integral are all equal. Therefore,
%----
$$ \Vert {\rm ad}A_1 ...  {\rm ad} A_n  U(t) \Vert
\leq  |t |^n  \prod _{j=1}^n
h    \varepsilon _{\psi(j)} (t).
 $$
%---
With  notations in Theorem 1.2 of \cite{A-L-N}, this amounts to the statement  of  Proposition \ref{p5.3}.

\hfill \carre

The Proposition below ends the proof of the first two claims in Theorem \ref{t1.1}.
We make use of the notation $\{ F, G \}$ for the  Poisson bracket of two functions  $F$ and $G$ on $H^2$ taking values
in ${\cal L} ( ( {\C} ^2 )^{\otimes N} ) $, the product in the Poisson bracket   being a matricial product.

\begin{prop}\label{p5.4}  i) For  $0<h<1$, the operator
$U_h^{red}(t)$ is  of the form  $Op_h^{weyl} ( U(t , \cdot , h)) $ where
 $U (t,  \cdot , h)$ is   matricial symbol (taking values
in ${\cal L} ( ( {\C}  ^2 )^{\otimes N} ) $) which is $C^{\infty}$
 on $H^2$ and which satisfies, for all multi-index $(\alpha , \beta)$,
 %----
\be\label{5.16}\Big | \partial_q ^{\alpha} \partial_p ^{\beta}  U ( t, q, p , h) \Big |
 \leq M (h, t) \prod _{mn} (   |t| \varepsilon _{mn} (t))^{|\alpha _{mn} + \beta _{mn} |}\ee
 %----
where $\varepsilon _{mn} (t)$ is the sequence of the Proposition \ref{p5.1} and
 %----
\be\label{5.17} M (h, t)  = \prod _{mn} ( 1+ K S^2 h |t|^2 \varepsilon  _{mn} (t)^2 )\ee
%---
with $K$ being a universal constant  and
%---
$$ S = \sup _{mn} \max (1 ,  |t| \varepsilon _{mn} (t) ). $$
%-----
In particular, the function $U (t,  \cdot , h)$ belongs to $S_{\infty}^{mat} ({\cal B}, |t| \varepsilon (t))$.

ii) Denoting  by $H _{int} ^{free} (t , q, p)$ the Weyl (matricial) symbol  of $H _{int} ^{free} (t)$ (defined in (\ref{5.1})), we have,
%---
\be\label{5.18}{d\over dt } U (t,  q, p  , h) = -i H _{int} ^{free} (t , q, p) U (t,  q, p , h)
- {h\over 2}  \{ H _{int} ^{free} (t , \cdot ) , U (t,  \cdot  , h)\}  (q , p).\ee
%-------
\end{prop}
The constant $K$ in (\ref{5.17}) is related to the Weyl calculus but not to the physical system.

{\it Proof of i)}  One first applies Theorem 1.2 of \cite{A-L-N} with
$A_h $ being components of the matricial operator  $ U_h^{red} (t)$ and $M=1$,
$m= 2$. The basis denoted $(e_j)$ in \cite{A-L-N} here is the basis ${\cal B} =(f_{mn})$ and the sequence
$(\varepsilon _j)$ is $|t| \varepsilon_{mn} (t)$. From   Proposition \ref{p5.3}, the hypotheses
(14)  of this Theorem is fulfilled.  From  Theorem 1.2 of \cite{A-L-N},
there exists a matricial function  $U ( t, \cdot , h)$ whose components lies in
$S_2 ( {\cal B} , |t| \varepsilon (t))$ such that $U_h^{red} (t) = Op_h^{weyl} ( U ( t, \cdot , h) )$.
Moreover, this function satisfies (\ref{5.16}) for $\alpha = \beta = 0$.
Next, for each multi-index $(\alpha , \beta)$, one applies
Theorem 1.2 of \cite{A-L-N} to the components of  the following operator and constant $M$,
%----
$$ A_{\alpha \beta} =
 ({\rm ad} \widetilde P_h )^{\alpha} ({\rm ad} \widetilde Q_h )^{\beta}  U_{h}^{red}  (t),\qquad M = \prod _{mn} (  h |t| \varepsilon _{mn} (t))^{|\alpha _{mn} + \beta _{mn} |} $$
%--------
with $\widetilde P_h  ( f_{mn} ) = P_h  ( f_{mn} ) \otimes I$. According to  Theorem 1.2 of \cite{A-L-N},
there exists  a matricial function  $U _{\alpha \beta} ( t, \cdot , h)$ with components belonging to
$S_2 ( {\cal B} , |t| \varepsilon (t))$ such that $A_{\alpha \beta}  = Op_h^{weyl} ( U _{\alpha \beta} ( t, \cdot , h) )$.
Moreover, this function satisfies,
%----
$$ | U _{\alpha \beta} ( t, \cdot , h) | \leq  M (h, t) \prod _{mn} (  h |t| \varepsilon _{mn} (t))
^{|\alpha _{mn} + \beta _{mn} |}. $$
%-----
Using symbols for the composition formulas that may be applied when one of the operators in the composition has
its symbol defined as a continuous linear form on $H^2$ (see
Proposition 2.6 of \cite{A-L-N}),
%----
$$  U_{\alpha \beta  } ( t, \cdot , h) = c_{\alpha \beta} h^{|\alpha + \beta |}
\partial_q ^{\alpha} \partial_p ^{\beta}  U ( t, \cdot , h) $$
%---
with $|c_{\alpha \beta } | = 1$. The proof is completed.

{\it Proof of ii)}  According to point i) with Proposition 3.2 and  Proposition 2.8 of \cite{A-L-N}, we can extend equality (\ref{5.11}) for $f$ belonging to $W_1\otimes  ( {\C} ^2  )^{\otimes N} $. Consequently, for
all $f$ and $g$ in ${\cal D} \otimes  ( {\C} ^2  )^{\otimes N} $,
%-----
$$ {d\over dt } < U_h^{red} (t)  f , g> = i < U_h^{red} (t) f , H_{int}^{free} (t) g>. $$
%-----
The  right hand side $\Phi_t (q , p, h)$ of (\ref{5.18}) satisfies hypotheses $(H_1)$ and $(H_2)$ of
\cite{A-L-N}. From Theorem 2.2 of \cite{A-L-N}, one then may  associate
a quadratic form   $Q_h^{weyl } (\Phi_t)$ on ${\cal D}$ with it.  In view of Proposition 2.6 in \cite{A-L-N},
%-----
$$ - i < U_h^{red} (t) f , H_{int}^{free} (t) g>  = Q_h^{weyl } (\Phi_t) (f , g).$$
%-----
Point ii) then follows.

\hfill$\Box$

Thus, the first two claims in Theorem \ref{t1.1} are derived. We next study the last claim.

For each $t\in \R$, the nonnegative quadratic form $g_t $ on $H^2$ is defined setting,
%---
\be\label{g} g_t (X) = 3N  |t| \sum _{m=1 }^3 \sum _{\lambda =1}^N \int_0^t | X \cdot
B_{m, a_{\lambda }, s }  | ^2 ds, \ee
%-----
where the  $B_{m, a_{\lambda }, s } $ are elements of $H^2$ defined in (\ref{5.3})
when identifying $H^2$ and $H_{\C}$.

\begin{prop}\label{p5.5} The matricial Weyl symbol $U(t, q, p, h)$ of the operator $U_h^{red}(t)$
is $C^{\infty}$ on $H^2$ and satisfies,
for all integers $m$, for all vectors $ X_1 , \dots, X_m $ in $H^2$, for all $t\in \R$,
%----
\be\label{5.19} | d^m U (t, q, p, h) ( X_1 , \dots, X_m ) | \leq M(h, t)\prod _{j=1}^m g_t (X_j ) ^{1/2}
\leq  M(h, t)(\prod _{j=1}^m K t |X_j|)\ee
%----
where $M(h, t)$ is defined in (\ref{5.17}) and $K$ is a constant.
\end{prop}

{\it Proof.} For any  $t\in \R$ and for each $X\in H^2$, set
%----
$$ N_t (X) =  \sum _{m=1 }^3 \sum _{\lambda =1}^N
 |  X \cdot B_{m, a_{\lambda }, t }  |. $$
 %-----
Set, for any $X$ and $(q, p)$ in $H^2$, $F_X (q , p)= \sigma ( (q, p), X)$,
 where $\sigma $ is the symplectic form.
For all $X$ in $H^2$ and all
$t\in \R$, one has
%---
$$ \Vert [ Op_h^{weyl} (F_X ) , H_{int}^{free} (t) ]  \Vert \leq h  N_t (X). $$
%---
Equality (\ref{5.14}) together with the above estimate show that, for all vectors $Y_1,\dots,Y_p$
in $H^2 $,
setting
 $A_{Y_j} = Op_h^{weyl} (F_{Y_j} ) \otimes I$,
%--
$$\Vert {\rm ad} A_{Y_1} \dots {\rm ad} A_{Y_p}
U_h^{red} (t) \Vert \leq h^p  \sum _{\varphi \in S_p}  \int _{  \Delta _ p(t)}
N_{s_1} (  Y_{\varphi (1)}) \cdots N_{s_p} (  Y_{\varphi (p)})
 ds_1 \cdots ds_p$$
%----
\be\label{5.XXX} \leq h^p  \prod _{j=1}^p  \int _0^t N_{s} (  Y_{j}) ds
 \leq   h^p  \prod _{j=1}^p  g_t (Y_j) ^{1/2}.\ee
%------

Fix $X_1,\dots, X_m $ and let $(\alpha , \beta )$ be a multi-index. Let
$A_{\alpha \beta }$ be defined in the proof of Proposition \ref{p5.4}.
Note that, for all elements $f_{mn}$ in the basis ${\cal B}$, one has,
%---
$$ \sup _{ |s| < |t| }  N_s ( (f _{mn }, 0 ) \leq \varepsilon _{mn} (t),\qquad \sup _{ |s| < |t| }  N_s ( (0, f _{mn } ) \leq \varepsilon _{mn} (t).   $$
%-----
The above inequality shows that,
%----
$$ \Vert A_{\alpha \beta } {\rm ad}  Op_h^{weyl} (F_{X_1} ) \otimes I\cdots {\rm ad}  Op_h^{weyl} (F_{X_m} )\otimes I
U_h^{red} (t) \Vert \leq  \prod _{mn }
( h |t | \varepsilon _{mn} (t) )^{\alpha _{mn} + \beta _{mn} }
\ \prod _{j=1}^m \int _0^t h N_{s} (  X_{j}) ds. $$
%---
Set
%---
$$ V(t) = ({\rm ad} Op_h^{weyl} (F_{X_1} ) \otimes I ) \cdots ({\rm ad} Op_h^{weyl} (F_{X_m} ) \otimes I )
 U_h^{red} (t).$$
 %----
The characterization theorem, Theorem 1.2 of \cite{A-L-N}, shows that $V(t)$ is associated with a function $V(t, q, p)$ by the Weyl calculus. In particular, it satisfies,
%----
$$ | V(t, q, p)| \leq M(h, t) \prod _{j=1}^m \ \int _0^t h N_{s} (  X_{j}) ds,  $$
%---
where  $M(h, t)$ is defined by (\ref{5.17}). Using composition  formulas that are valid if one of the operators is a linear form, one knows that,
%---
$$ V(t, q, p) =  h^m i^{-m} d^m U (t, q, p, h) ( X_1 , \dots, X_m ). $$
%-----
Consequently,
%----
$$| d^m U (t, q, p, h) ( X_1 , \dots, X_m ) |  \leq M(h, t) \prod _{j=1}^m \
 \int _0^t  N_{s} (  X_{j}) ds  \leq M(h, t) \prod _{j=1}^m  g_t (X_j) ^{1/2}. $$
%------
The first inequality in (\ref{5.19}) then follows.

Since $g_t (X) ^{1/2} \leq K |t| |X|$ then the second inequality in (\ref{5.19}) is valid.

\hfill$\Box$

{\it End of the Proof of Theorem \ref{t1.1} $(iv)$.} Extending the differential of order $m$ to a multilinear form on
the complexified  $H_{\C}^2$, the extension may be  defined by,
%---
$$  {\cal U} (t, X, h)  = \sum _{j\geq 0} {1\over j!} d^j U (  t, 0, 0, h)
( X )^j, \qquad X = (q , p)\in H_{\C}^2.$$
%---
The estimate (\ref{1.2}) comes from the fact that one may also write,
%---
$$  {\cal U} (t, X+ i Y, h)  = \sum _{j\geq 0} {1\over j!} d^j U (  t, X, h)
( iY )^j, $$
%---
for $X$ and $Y$in $H^2$.

\section{Proof of Theorem \ref{t1.1} $(iii)$
(Semiclassical expansion of the reduced propagator symbol).}\label{s6}

 We study here the Weyl symbol $U(t, q, p, h)$ ($t \in \R$, $(p, q) \in H^2$)
of the operator $U^{red}_h (t)$ of (\ref{1.1}) or (\ref{5.6}). We aim to derive an asymptotic expansion of this symbol as $h$ goes to $0$,
$$
 U(t, q, p, h) \sim \sum _j g_j(t, q, p)h^j .
$$

 The functions $g_j(t, q, p)$ takes values in ${\cal L} ( ({\C}^2) ^{\otimes N})$ similarly as the symbol itself. Since
the symbol  $U (t, q, p, h)$ satisfies (\ref{5.18}), it is natural that the functions $g_j$ satisfies,
%-----
\be\label{6.1a}{d\over dt} g_0(t, q, p) = - i H^{free}_{int}(t, q,p) g_0 (t, q, p),\qquad g_0 (0, q, p) = I,
\ee
%------
where  $H^{free}_{int}(t, q,p)$ is the matricial Weyl symbol of  $H^{free}_{int}(t)$ defined in (\ref{5.1}). For all  $j\geq 1$, one should have,
%-----
\be\label{6.1b}{d\over dt} g_j(t, q, p) = - i H^{free}_{int}(t, q,p) g_j (t, q, p)
- {1\over 2}  \{ H _{int} ^{free} (t , \cdot ) , g_{j-1} (t,  \cdot  , h)\}  (q , p)  ,
\qquad g_j  (0, q, p) = 0.\ee
%------

\begin{prop}\label{p6.1} For all $(q, p)$ in $H^2$, the solutions  $t \mapsto    g_j(t, q, p)$
 of the system (\ref{6.1a})(\ref{6.1b}) (taking values in ${\cal L} (( \C^2  )^{\otimes N} )$) are well defined on $\R$,
and  $ g_0(t, q, p)$ is unitary. Moreover, one has
for all multi-indices $(\alpha , \beta)$,
%----
\be\label{6.2} \Big | \partial _q ^{\alpha } \partial _p^{\beta} g_j(t, q, p) \Big |
\leq  \prod _{mn} (  K |t|\varepsilon  _{mn } (t) ) ^{ \alpha _{mn}
+ \beta _{mn}  },\ee
%------
 where $ \varepsilon _{mn } (t)$ is defined in  Proposition \ref{p5.1},  for some positive real number $K$. The function  $g_j(t, \cdot )$ then belongs to $S_{\infty }^{mat} ( {\cal B}, K |t| \varepsilon (t))$.
\end{prop}

{\it Proof.} The existence of time global solution  and the unitarity of $g_0(t, q, p)$
are standard. Let us prove the bounds on $g_0$. Let $  P_1  , \dots,  P_r$  be a finite sequence of operators where each operator $P_j$
is either one of the operators ${\partial \over \partial q_{mn}}  $ or  ${\partial \over \partial p_{mn}}  $,
the corresponding index  $(m, n)$ being noted $\psi (j)$. Applying the operator $P_1$
on the two sides of (\ref{6.1a}) and solving the resulting system, we see that,
%----
\be\label{6.3} P_1  g_0(t, q, p) = -i \int _0^t g_0(t, q, p) g_0(s, q, p)^{\star } (P_1H^{free}_{int}(s, q,p) ) g_0(s, q, p) ds.\ee
%----
Iterating, one obtains,
%----
$$P_1 \dots P_r g_0(t, q, p) = (-i)^r \sum _{\varphi \in S_r}  \int _{  \Delta _ r(t)}
g_0(t, q, p) g_0(s_r , q, p)^{\star } \Big (  P_{\varphi (r)} H^{free}_{int}(s_r, q, p) \Big ) $$
%---
\be\label{6.4} g_0(s_r, q, p) g_0(s_{r-1}  , q, p)^{\star } \Big (  P_{\varphi (r-1)} H^{free}_{int}(s_{r-1} , q, p) \Big )
\ee
%-----
$$ g_0(s_2, q, p) g_0(s_{1}  , q, p)^{\star } \Big (  P_{\varphi (1)} H^{free}_{int}(s_{1} , q, p) \Big )
g_0( s_1 , q, p) \, ds_1\dots ds_r $$
%----------
where
  $S_r$ is the set of all bijections $\varphi $  of
the set $\{ 1, ... , r \}$  and where  $ \Delta _ r(t)$ is defined in (\ref{5.15}), if $t>0$.
Since $g_0(t, q, p)$ is unitary, one observes using (\ref{5.2}),
%-----
$$ \Big | P_1 \dots P_r g_0(t, q, p) \Big | \leq \sum _{\varphi \in S_r}
\prod _{j=1}^r \varepsilon _{\psi ( \varphi (j))}(t)
 \int _{  \Delta _ r(t)} ds_1\cdots ds_r. $$
%----------
The factors in front of the integral are all equals. The value of the integral is
 $|t|^r / r!$. Consequently,
%--------
$$ \Big | P_1 \dots P_r g_0(t, q, p) \Big | \leq |t|^r
\prod _{j=1}^r \varepsilon _{\psi ( j)}(t). $$
%----
Turning back to usual  notations  for  multi-indices, one obtains
(\ref{6.2}) for $g_0$. Now suppose that the estimates (\ref{6.2}) are valid for the integer $j-1$. Equations (\ref{6.1a})(\ref{6.1b}) and Duhamel principle show that,
%-----
$$ g_j(t, q, p) = - {1\over 2}  \int_0^t g_0  (t , q, p)
g _0  (s , q, p)^{\star }
 \{ H _{int} ^{free} (t , \cdot ) , g_{j-1} (t,  \cdot  , h)\}  (q , p).
$$
%------
The iteration hypothesis gives, if $|s| < |t|$,  $g_j(s, \cdot)$ belongs to
 $S_{\infty } ^{mat} ( {\cal B},   |t| \varepsilon(t))$. According to Proposition \ref{p3.2},  the Poisson bracket
$ \{  H^{free}_{int}(t, \cdot )\, , \, U(t, \cdot , h )\}   (q , p)   $
belongs to a set $S_{\infty } ^{mat} ({\cal B} , K' |t| \varepsilon(t)) $.
The function $g_0(t, \cdot )$ has already been seen to be in $S_{\infty}^{mat} ( {\cal B},  |t| \varepsilon(t)) $.
The estimates on $g_j$ then follows.

\hfill $\Box$

\begin{prop}\label{p6.2}   The family of unitary matrices   $ g_0( t, q , p)$,
solution to (\ref{6.1a}) and taking values in ${\cal L} (( \C^2  )^{\otimes N} )$, satisfies, for all $s$ and $t$ in $\R$, for
any $(q, p)$  in $H^2$,
%----
\be\label{6.5} g_0(t, q, p)\  g_0(s , q, p)^{\star} =  g_0( t-s , \chi_s (q , p)).\ee
%--------
\end{prop}

{\it Proof.}  For any $s\in \R$, the two sides of  the equality
are solutions to the same differential system,
%----
$$ {d\over dt} F(t, q, p) = i H_{int}^{free} (t) F(t, q, p ),\qquad F(s, q, p) = I. $$
%--
For the right hand side, one uses the fact  that
$H_{int}^{free } ( t-s , \chi_s (q , p)) )=  H_{int}^{free } ( t, q, p) $.
The Proposition then follows.

 \hfill \carre

\begin{prop}\label{p6.3} Let $U(t, q, p, h)$ be the
symbol of $U^{red}_h (t)$ and $g_j(t, q, p)$ ($j\geq 1$) be the solutions to $(\ref{6.1a})(\ref{6.1b})$. Then
we have for all integer $m$,

\be\label{6.6} U(t, q, p, h ) - \sum _{j=0}^m  g_j (t, q, p)h^j  =
h^{m+1}   R_m(t, q, p, h),\ee
%------
where $R_m( t , \cdot , h)$ lies in a set  $S_{\infty } ^{mat} ({\cal B} , K |t| \varepsilon(t)) $ and with
$K$ being some constant.
\end{prop}

{\it Proof.} Set,
%---
\be\label{s}S_m (t, q, p, h) =  \sum _{j=0}^m  g_j (t, q, p)h^j.\ee
%---
For $m= -1$, we agree that $S_{-1} (t, \cdot , h)= 0$ and $R_{-1} (t, \cdot , h)= U(t, \cdot , h)$. The bounds on $R_{-1}$ are those of Theorem \ref{t1.1}. Suppose that the Proposition is proved for the integer $m-1$. With notations (\ref{6.6}) and (\ref{s}), equation (\ref{5.18}) reads as,
%----
$$ {d\over dt } U (t,  \cdot  , h) = -i H _{int} ^{free} (t , \cdot) U (t,  \cdot , h)
- {h\over 2}  \{ H _{int} ^{free} (t , \cdot ) , S_{m-1} (t,  \cdot  , h) + h^m
R_{m-1} (t,  \cdot  , h) \}.   $$
%-------
Since the $g _j (t , q, p)$ satisfy (\ref{6.1a}) and (\ref{6.1b}),  one has,
%----
$$ {d \over dt }   S_m (t, \cdot, h)  =
-i H _{int} ^{free} (t , \cdot )   S_m (t, \cdot , h) - {h\over 2}  \{ H _{int} ^{free} (t , \cdot ) , S_{m-1} (t,  \cdot  , h)\}.  $$
%---
Since $g_0$ satisfies (\ref{6.5}) and since $ U- S_m $ is vanishing at $t= 0$,
 Duhamel principle gives that,
%------------------
$$ U(t, q, p, h)  - S_m (t, q , p , h)    = - {h^{m+1}  \over 2} \int_0^t g_0  (t , q, p)
g _0  (s , q, p)^{\star }
\{   H^{free}_{int}(s, \cdot )\, , \, R_{m-1} (s, \cdot , h)\}  (q , p)  \, ds.$$
%------------------------------
According to iteration hypothesis, the function $R_{m-1} (s, \cdot , h) $ belongs to
 $S_{\infty } ^{mat} ( {\cal B},   |t| \varepsilon(t))$. In view of Proposition \ref{p3.2},  the Poisson bracket
$ \{  H^{free}_{int}(t, \cdot )\, , \, R_{m-1} (s, \cdot , h) \}   (q , p)   $
belongs to a set  $S_{\infty } ^{mat} ({\cal B} , K' |t| \varepsilon(t)) $.
According to Proposition \ref{p6.1}, the function $g_0(t, \cdot )$ est dans $S_{\infty}^{mat} ( {\cal B},  |t| \varepsilon(t)) $.
The proof of the Proposition is completed

\hfill $\Box$

\section{Proof of Theorem \ref{t1.2}. (Consequences of  analyticity properties.)}\label{s7}

\begin{theo}\label{t7.1} Let $F$ be a function in the set $S_{\infty}  ^{mat}({\cal B},  \varepsilon)$.
Suppose that,  for all integers $m$, for any vectors $ X_1 , .. X_m $ in $H^2$,
%----
\be\label{7.1}  | d^m F ( q, p) ( X_1 , .. X_m ) | \leq M \prod _{j=1}^m (  K |X_j|).\ee
%----
Then, we have for all $X$ and $Y$ in $H^2$, for any
$a$ and $b$ in $({\C}^2)^{\otimes N}$ with norm $1$,
%---
\be\label{7.2}\left |  < Op_h^{weyl} (F) (\Psi_{X , h} \otimes a ), (\Psi_{Y  , h} \otimes b )>  \right |
 \leq M  e^{  K |X-Y|- {1\over 4h} |X-Y|^2 }.\ee
%-------
\end{theo}

{\it Proof.} Under our hypotheses, the function $F$  has a stochastic extension
$\widetilde F$  in $L^1 (B^2 , \mu _{B^2, t}  )$ (see Section \ref{s2.B}), for all $t>0$. For any $t>0$, the operator,
%---------
\be\label{7.3} (H_tF)  (X) =
 \int _{B^2}  \widetilde F( X +Y ) d\mu _{B^2, t} (Y)\ee
%----
is considered as the heat operator.
The function $H_{h/2} F(  \cdot )$ also satisfies,
%----
$$ | d^m H_{h/2} F  ( q, p) ( X_1 , \dots,  X_m ) | \leq M \prod _{j=1}^m (  K |X_j|).  $$
%----
Consequently, as in the proof of Theorem \ref{t1.1} (end of Section \ref{s5}),
 functions $F$ and $H_{h/2} F$ are extended  to
holomorphic functions  on $H_{\C}^2$, taking values in ${\cal L} ( ( {\C}^2) ^{ \otimes N} )$, denoted
${\cal F}$ and  $H_{h/2}  {\cal F} ( \cdot, h)$ and satisfying,
%----
\be\label{7.4} |  H_{h/2}  {\cal F} ( X+ i Y , h) |  \leq M  e^{  K  |Y |},\qquad X+iY \in H_{\C}^2.\ee
%---
The bi-symbol of $Op_h^{weyl} (F)$ defined in (\ref{3.4})  is,
%---
\be\label{7.5} S_h ( Op_h^{weyl} (F) ) ( q , p, q' , p') = \Big ( H_{h/2} {\cal F} \Big )
 \left ( { q+ip \over 2} +  { q' -ip' \over 2} ,  { p-iq \over 2} +  { p'+iq' \over 2}, h
 \right ).\ee
 %----
 Indeed, one knows that the bi-symbol of an operator is a function on $H^2 \times H^2$ which is holomorphic
 with respect to the first variable and anti-holomorphic with respect to the seconde variable $H^2$ identified with $H_{\C}$. The above function  is the only one sharing theses properties and
whose  restriction to the diagonal  equals to $H_{h/2} F(  \cdot )$, which is the  Wick symbol
  of the operator $Op_h^{weyl} (F)$ (from \cite{A-J-N}, Theorem 7.1).
For all $X$ and $Y$ in $H^ 2$, for any  $a$ and $b$ in $E$ with unit norms,  one deduces
from (\ref{3.4}),
 %----
 $$   < Op_h^{weyl} (F) (\Psi_{X , h} \otimes a ) ,  (\Psi_{Y  , h} \otimes b )>
 = S_h ( Op_h^{weyl} (F) ) ( X, Y )(a) \cdot b \ < \Psi_{X , h} , \Psi_{Y , h}>, $$
 %---
thus, setting $X = (q, p)$, $Y = (q', p')$,
 %------
 $$ = \Big ( H_{h/2} {\cal F} \Big )
 \left ( { q+ip \over 2} +  { q' -ip' \over 2} ,  { p-iq \over 2} +  { p'+iq' \over 2}, h
 \right ) (a) \cdot b \   \ < \Psi_{X , h} , \Psi_{Y , h}>. $$
 %----
From (\ref{7.4}) and using a standard  equality concerning the scalar product of
 coherent states (\cite{A-L-N}, formula (2.3)), we have indeed (\ref{7.2}).

\hfill $\Box$

{\it Proof of Theorem \ref{t1.2}.} One has,
%----
$$  < e^{-i{t \over h}  H(h)} (\Psi_{X , h} \otimes a ), (\Psi_{Z  , h} \otimes b )>  =
< e^{-i{t \over h} H_0} U_h^{red} (t) (\Psi_{X , h} \otimes a ), (\Psi_{Z  , h} \otimes b )> $$
%---
$$ = <  U_h^{red} (t) (\Psi_{X , h} \otimes a ), (\Psi_{\chi_{-t} (Z)  , h} \otimes b )>. $$
%---
It is then sufficient to apply   Theorem  \ref{t7.1} with $F (q , p) = U(t, q, p, h)$.

\hfill $\Box$

\medskip

laurent.amour@univ-reims.fr\newline
LMR EA 4535 and FR CNRS 3399, Universit\'e de Reims Champagne-Ardenne,
 Moulin de la Housse, BP 1039,
 51687 REIMS Cedex 2, France.

rlascar@math.univ-paris-diderot.fr\newline
Institut Mathématique de Jussieu UMR CNRS 7586,  Analyse Algébrique, 4 Place Jussieu, 75005 Paris, France.

jean.nourrigat@univ-reims.fr\newline
LMR EA 4535 and FR CNRS 3399, Universit\'e de Reims Champagne-Ardenne,
 Moulin de la Housse, BP 1039,
 51687 REIMS Cedex 2, France.

\end{document}